\documentclass[reqno,12pt]{amsart}

\newtheorem{theorem}{Theorem}[section]
\newtheorem{lemma}[theorem]{Lemma}
\newtheorem{corollary}[theorem]{Corollary}

\theoremstyle{definition}

\theoremstyle{remark}
\newtheorem{remark}[theorem]{Remark}

\makeatletter
\def\dashint{\operatorname%
{\,\,\text{\bf--}\kern-.98em\DOTSI\intop\ilimits@\!\!}}
\makeatother

\newcommand\bA{\mathbb{A}}
\newcommand\bB{\mathbb{B}}
\newcommand\bH{\mathbb{H}}

\newcommand\bR{\mathbb{R}}
\newcommand\bQ{\mathbb{Q}}
 
 \newcommand\cA{\mathcal{A}}
\newcommand\cB{\mathcal{B}}

\newcommand\cH{\mathcal{H}}
 
\newcommand\cL{\mathcal{L}}

\newcommand{\Div}{{\rm div}\,}
\newcommand{\osc}{{\rm osc}\,}

\newcommand{\cHO}{\overset{\scriptscriptstyle0}%
{\mathcal H}\,\!}
\newcommand{\WO}{\overset{\scriptscriptstyle0}%
{W}\,\!}

 \newcommand{\mysection}[1]{\section{#1}
 \setcounter{equation}{0}}

\begin{document}

\title[Parabolic   equations with VMO coefficients]
{Parabolic   equations with VMO coefficients
in spaces with mixed norms}
\author{N.V. Krylov}
\thanks{The work  was partially supported by  
NSF Grant DMS-0140405}
\email{krylov@math.umn.edu}
\address{127 Vincent Hall, University of Minnesota, 
Minneapolis, MN, 55455}
 
 

\begin{abstract}
An $L_{q}(L_{p})$-theory
 of divergence and non-divergence form
  parabolic
equations is presented.
The main coefficients are supposed to belong to the
class $VMO_{x}$,
which, in particular, contains all measurable functions
depending only on $t$. 
The method of proving simplifies the methods
previously used in the case $p=q$. 
\end{abstract}

\maketitle

\mysection{Introduction}
                                                      \label{intro}

The goal of this article is to
prove the solvability of parabolic second-order
divergence and non-divergence type equations
in Sobolev spaces with mixed norms.

More precisely, we are dealing with two types of parabolic operators:
$$
Lu(t,x)=u_{t}(t,x)+a^{ij}(t,x)u_{x^{i}x^{j}}(t,x)+
b^{i}(t,x)u_{x^{i}}(t,x)+c(t,x)u(t,x) ,
$$
$$
\cL u(t,x)=u_{t}(t,x)+\big(a^{ij}(t,x)u_{x^{i}}(t,x)+
\hat{b}^{j}(t,x)u(t,x)\big)_{x^{j}} 
$$
$$
+b^{i}(t,x)u_{x^{i}}(t,x)+c(t,x)u(t,x)
$$
acting on functions given on 
$$
\bR^{d+1}=\{(t,x):t\in\bR,x\in\bR^{d}\},
$$
where $\bR^{d}$ is a $d-$dimensional Euclidean space of points
$x=(x^{1},...,x^{d})$.

The interest in results concerning equations
in spaces with mixed $L_{q}(L_{p})$-norms
arises, in particular, when one wants 
to have better regularity of traces
of solutions for each $t$ while treating
linear or nonlinear equations
(see, for instance, \cite{GS} and \cite{MS}
for applications to the Navier-Stokes
equations).

Parabolic equations in $L_{q}(L_{p})$-spaces  
have been investigated in many articles
for at least forty years.
The interested reader can find many references
and discussions of methods and obtained results
in \cite{DHP}, \cite{HHH}, and \cite{Kr02}.

However, it seems to the author that apart from
\cite{Kr02} (also see the references therein)
in most other papers concerning $L_{q}(L_{p})$-spaces
the methods heavily depend on the properties
of the elliptic part in $L$ or $\cL$,
which is supposed to be independent of
$t$ and have well behaving resolvent
or generate a ``good" semigroup.
However, in \cite{Am} (also see references therein)
 there is a general
theorem allowing one to treat the case when
the coefficients are continuous in $t$.
These restrictions exclude
parabolic equations with coefficients measurable or even
VMO in $t$
(even if they are independent of $x$,
the case considered in \cite{Kr02}).
In particular, in \cite{HHH} the authors only consider
equations with VMO coefficients independent
of time, although combining
their results with \cite{Am} would include equations with
coefficients continuous in $t$. 
By  the way, in the particular case
that $q=p$ this also does not allow one to cover
the results of \cite{BC}, where
the coefficients  are in $VMO(\bR^{d+1})$.
Speaking about the case $q=p$ it is worth saying
that there is a quite extensive literature
about equations and systems with VMO coefficients.
The interested reader can consult \cite{By1}, \cite{By2},
 \cite{DHP}, \cite{Gu}, \cite{HHH}, \cite{MPS}, \cite{PS},
\cite{PS1}, \cite{PS2}, \cite{RT},
\cite{So1}, \cite{So2}, \cite{So3},
and the references therein.

Our approach is based on a method
from \cite{Kr06} and further developed in
\cite{KDK1} and \cite{KDK2}, where everything
hinges on   a priori {\em pointwise\/} estimates
of the sharp functions of the second-order 
spatial derivatives of solutions.
This method allows one to avoid
using generalizations of the Calder\'on-Zygmung
theorem and the 
Coifman-Rochberg-Weis commutator theorem
as is often done when VMO is involved
(see, for instance, \cite{DHP}, \cite{Gu}, \cite{HHH}, \cite{MPS},
\cite{PS},
\cite{PS1},\cite{PS2}, \cite{RT},
\cite{So1}, \cite{So2}, \cite{So3},
and the references therein).
However, it is worth noting that if $p=q$
 there is
an approach to the {\em divergence\/} type equation
 suggested in \cite{By1} and \cite{By2},
which also does not use the above mentioned tools.
The approach from the present article has been already used
in a very interesting article \cite{Ki1}
to prove the solvability
in usual Sobolev spaces of parabolic equations 
with partially VMO coefficients when most of the coefficients
are just measurable in time and one of space variables
and VMO with respect to the others.

In \cite{Kr06}, in each small cylinder,
 the solution is split into two parts:
a function, that is ``harmonic" with respect
to the operator with ``frozen" coefficients,
 and the remainder.
In order to do this decomposition one has to know
that the corresponding boundary-value problems are solvable.
This is not very convenient if one has in mind higher-order
equations. 

It turns out that instead one can use
splitting of the right-hand side of the equation
and rely on solvability of equations in the whole space.
This approach not only simplifies some proofs from
\cite{Kr06} but also allows one to 
make stronger main technical estimates
(see Lemmas \ref{lemma 06.7.26.2} and \ref{lemma 06.10.9.5}),
which after being combined with an approach suggested
in \cite{Kr02} leads to $L_{q}(L_{p})$-theory.
Although, we are dealing only with the Cauchy problem
for second-order operators,
it seems   that the new technique, 
which we develop here, is applicable to higher-order
equations, systems, and boundary-value problems
for elliptic and parabolic equations with VMO
coefficients.

We are assuming that the main coefficients are
measurable in time and $VMO$ in spatial variables
and prove the solvability in $L_{q}(L_{p})$
spaces for $L$ if $q\geq p$ (Theorem \ref{theorem 7.20.3})
 and for $\cL$
without this restriction (Theorem \ref{theorem 7.20.4}).
Theorem \ref{theorem 7.20.3} generalizes the 
{\em corresponding\/} result of \cite{HHH}
to cover time-dependent coefficients.
However, note that the results in \cite{HHH}
are proved also for higher-order parabolic systems,
 arbitrary $p,q\in(1,\infty)$,
and  $L_{p}$-spaces with $A_{p}$ Muckenhoupt weights.

The paper is organized as follows.
In Section \ref{section 7.22.1} we state our main results.
Theorem \ref{theorem 7.20.3} and \ref{theorem 7.20.4}
are proved in Sections \ref{section 06.10.22.1} and
\ref{section 06.10.22.2}, respectively,
on the basis of Lemmas \ref{lemma 06.7.26.2}
and \ref{lemma 06.10.9.5}, respectively,
which are proved later. In Sections
\ref{section 7.22.2} and \ref{section 7.22.20}
we present our new approach to treating parabolic
equations with $VMO_{x}$ coefficients.
The main results of these two sections are
Theorems \ref{theorem 06.6.13.1} (non-divergence
equations) and \ref{theorem 06.10.6.1} (divergence equations)
about equations in usual Sobolev spaces without mixed norms.
If one takes functions independent of $t$,
these two theorems yield the basic estimates for
elliptic equations.
Finally, in Sections \ref{section 06.10.19.1}
and \ref{section 06.10.20.1} we prove Lemmas
\ref{lemma 06.7.26.2}
and \ref{lemma 06.10.9.5}, respectively.

 The work on this article was
stimulated by discussions during the author's
stay at 
 Centro di Ricerca
Matematica Ennio De Giorgi, Scuola Normale Superiore di Pisa,
and it is a great pleasure to bring my
sincere gratitude to G. Da Prato and M. Giaquinta
for the invitation and hospitality.

We finish the section introducing some notation.
Note that we also use without mentioning some common notation
from PDEs.
For 
$$
-\infty\leq S<T\leq\infty,\quad1<p, q<\infty
$$
 we set
$$
\bR_{S}=(S,\infty),\quad\bR^{d+1}_{S}=\bR_{S}\times\bR^{d},\quad
L_{p}=L_{p}(\bR^{d+1}),
$$
$$
L_{q,p}((S,T))=L_{q}((S,T),L_{p}(\bR^{d})),
\quad L_{q,p}=L_{q,p}(\bR),
$$
$$
W^{1,2}_{q,p}((S,T))=\{u:u,u_{t},u_{x},u_{xx}
\in L_{q,p}((S,T))\},
$$
$$
W^{1,2}_{q,p}=W^{1,2}_{q,p}(\bR),\quad
W^{1,2}_{p}(\bR^{d+1}_{S})=W^{1,2}_{p,p}(\bR_{S}),
\quad W^{1,2}_{p}=W^{1,2}_{p}(\bR^{d+1}),
$$
$$
\|u\|_{L_{q,p}((S,T))}^{q}
=\int_{S}^{T}\big(\int_{\bR^{d}}|u(t,x)|^{p}\,dx\big)^{q/p}\,dt,
$$
$$
\|u\|_{W^{1,2}_{q,p}((S,T))}=
\|u\|_{L_{q,p}((S,T))}+\|u_{x}\|_{L_{q,p}((a,b))}
$$
$$
+\|u_{xx}\|_{L_{q,p}((S,T))}
+\|u_{t}\|_{L_{q,p}((S,T))}.
$$
By $\WO^{1,2}_{q,p}((S,T))$ we mean the subspace
of $W^{1,2}_{q,p}(\bR_{S})$ consisting of functions $u(t,x)$
vanishing for $t>T$.  Finally,
$$
\cH^{1}_{q,p}((S,T))=(1-\Delta)^{1/2}W^{1,2}_{q,p}((S,T)),\quad
\cH^{1}_{q,p} =\cH^{1}_{q,p}(\bR),
$$
$$
\cHO^{1}_{q,p}((S,T))=(1-\Delta)^{1/2}\WO^{1,2}_{q,p}((S,T)),
$$
where $\Delta$ is the Laplacian in $x$ variables.
In the above notation we write $p$ in place of $q,p$ if $q=p$.
For instance, $\cHO^{1}_{p}((S,T))=\cHO^{1}_{p,p}((S,T))$.
In particular, $W^{1,2}_{p}(\bR^{d+1}_{S})=W^{1,2}_{p}(\bR_{S})$.
Finally,
$$
\bH^{-1}_{p}((S,T))=(1-\Delta)^{1/2}L_{p}((S,T)\times\bR^{d}),
 \quad\bH^{-1}_{p}=\bH^{-1}_{p}(\bR).
$$

\mysection{Main results}
                                      \label{section 7.22.1}

We assume that the coefficients of 
$L$ and $\cL$ are measurable
and by magnitude are dominated by a constant $K<\infty$. We also assume
that the matrices $a=(a^{ij})$ are, perhaps, nonsymmetric and satisfy
\begin{equation}
                                                   \label{7.18.1}
  a^{ij} \lambda^{i}\lambda^{j}\geq\kappa|\lambda|^{2}
\end{equation}
for all $\lambda\in\bR^{d}$ and all possible values of arguments.
Here $\kappa>0$ is a fixed constant.

To state our main assumption we set $B_{r}(x)$ to be the open ball 
in $\bR^{d}$ of radius $r$ centered at $x$, $B_{r}=B_{r}(0)$,
$Q_{r}(t,x)=(t,t+r^{2})\times B_{r}(x)$, $Q_{r}=Q_{r}(0,0)$,
$\bB$ the collection of open balls in $\bR^{d}$,
and $\bQ$ the collection of $Q_{r}(t,x)$, $(t,x)\in\bR^{d+1}$,
$r\in(0,\infty)$.
Denote
$$
\osc_{x}(a,Q_{r}(t,x))=
r^{-2}|B_{r}|^{-2}
\int_{t}^{t+r^{2}}\int_{y,z\in B_{r}(x)}|a(s,y)-a(s,z)|\,dydzds,
$$
$$
a^{\#(x)}_{R}=\sup_{(t,x)\in\bR^{d+1}}\sup_{r\leq  R}
\osc_{x}(a,Q_{r}(t,x)) .
$$
 We assume that $a\in VMO_{x}$, that is
\begin{equation}
                                                   \label{9.15.1}
\lim_{R\to0}a^{\#(x)}_{R}=0.
\end{equation}
For convenience of stating our results we take any 
increasing
continuous function $\omega(R)$ on $[0,\infty)$, such that
$\omega(0)=0$ and $a^{\#(x)}_{R}\leq\omega(R)$ for all
$R\in(0,\infty)$. 
Obviously, $a\in VMO_{x}$ if $a$ depends only on $t$.

Needless to say   all equations below are
understood in the sense of generalized functions.

Now we fix $T\in(0,\infty)$ and $q,p\in(1,\infty)$, set
$$
\Omega(T)=(0,T)\times\bR^{d}
$$
and state our main results.
\begin{theorem}
                                              \label{theorem 7.20.3}
Let  $q\geq p$. Then
for any $f\in L_{q,p}((0,T))$ there exists
a unique $u\in\WO^{1,2}_{q,p}((0,T))$ such that
$Lu=f$ in $\Omega(T)$. Furthermore, there is a constant $N$,
depending only on $d$, $T$, $K$, $\kappa$, $q,p$, and 
the function $\omega$, such that for any $u\in\WO^{1,2}_{q,p}((0,T))$
we have
\begin{equation}
                                               \label{7.21.6}
\|u\|_{W^{1,2}_{q,p}((0,T))}\leq N\|Lu\|_{L_{q,p}((0,T))}.
\end{equation}

\end{theorem}
 
\begin{remark}
                                      \label{remark 7.25.2}
Theorem \ref{theorem 7.20.3} is
similar to some results from \cite{Kr02} and \cite{HHH}
(also see the references therein).
However, in both articles there is no restriction on $p,q$.
On the other hand,
  in \cite{Kr02} the coefficients are independent
of $x$ and in \cite{HHH} they are independent of $t$.
As we have already pointed out in the Introduction,
by relying on \cite{Am},
some results from \cite{HHH} can be extended
to cover the case of coefficients continuous in $t$.
\end{remark}

\begin{theorem}
                                      \label{theorem 7.20.4}
Let $f=(f^{1},...,f^{d})$, $g,f^{i}\in L_{q,p}((0,T))$
for $i=1,...,d$. Then there is a unique
$u\in\cHO^{1}_{q,p}((0,T))$ such that 
$\cL u=\Div f+g$ in $\Omega(T)$.
Furthermore, there is a constant $N$,
depending only on $d$, $T$, $K$, $\kappa$, $q,p$, and 
the function $\omega$, such that 
\begin{equation}
                                            \label{7.25.1}
\|u \|_{L_{q,p}((0,T))}+
\|u_{x}\|_{L_{q,p}((0,T))}\leq N
\big(\|f\|_{L_{q,p}((0,T))}+\|g\|_{L_{q,p}((0,T))}\big).
\end{equation}
\end{theorem}

\begin{remark}
                                           \label{remark 7.23.1}
As usual in such situations,
from our proofs one can see that instead of the assumption
that
$a\in VMO_{x}$ we are, actually, using that
there   exists an $R\in(0,\infty)$ such that
$a^{\#(x)}_{R}\leq\varepsilon$, where $\varepsilon>0$
is a constant depending only on $d,p,\kappa,K$.
\end{remark}
\begin{remark}
                                                \label{remark 7.21.2}

Denote  
$$
 u _{Q }=\dashint_{Q } u(s,y)\,dyds,\quad
$$
the average value of a function $u(s,y)$ over $Q\in\bQ $ and
$$
u_{B }(t)=\dashint_{B } u(t,y)\,dy 
$$
the average value of a function $u(t,y)$ over $B\in\bB $.

Also introduce $\bA$ as the set of $d\times d$ matrix-valued
measurable
functions $a=a(t)$ depending only on $t$,
satisfying conditions \eqref{7.18.1} and such that
$|a^{ij}|\leq K$.

 A standard
fact to remember is
that for any $\bar{a}\in\bA$
$$
\osc_{x}(a,Q_{r})
\leq2 \dashint_{Q_{r}}|a(s,x)-\bar{a}(s)|\,dxds
$$
and for $\bar{a}(t)= a _{B_{r}}(t)$
$$
\dashint_{Q_{r}}|a(s,x)-\bar{a}(s)|\,dxds\leq\osc_{x}(a,Q_{r}).
$$
This allows one to give obvious equivalent definitions of
$VMO_{x}$.
\end{remark}

\mysection{Proof of Theorem \protect\ref{theorem 7.20.3}}
                                      \label{section 06.10.22.1}

 The following fact, which we prove in
Section \ref{section 06.10.19.1}, is a considerable improvement
of the key inequality from the proof
of Theorem 3.6 of~\cite{Kr06}.
It goes without saying that the assumptions
under which Theorem  \ref{theorem 7.20.3}
is stated are supposed to hold.

\begin{lemma}
                                             \label{lemma 06.7.26.2}
Let  $b=0$  and $c=0$. Then 
  there exists a constant $N=N(d,\kappa,K,p, \omega)$ such that
for any $u\in C^{\infty}_{0}(\bR^{d+1})$, $\nu\geq16$, and
$r\in(0,1/\nu]$     we have
\begin{equation}
                                                \label{06.7.26.8}
 (|u_{xx}-(u_{xx})_{Q_{r} }|^{p})_{Q_{r} }
\leq N \nu ^{d+2} \cA_{\nu r}+
  N(\nu^{-p}+\nu ^{d+2} 
  \hat{a}^{1/2})\cB_{\nu r},
\end{equation}
where
$$
\cA_{\rho}= (
|f|^{p})_{Q_{\rho} },
\quad
\cB_{\rho}=( |u_{xx}|^{p} )_{Q_{\rho} },
\quad
\hat{a}= a^{\#(x)}_{  \nu r },\quad f=Lu .
$$
 \end{lemma}

\begin{corollary}
                                         \label{corollary 06.7.28.1}
Let $b=0$ and $c=0$. Then
there exists a constant $N$ depending only
on $d,p,\kappa,K$, and $\omega $,
such that for any $u\in C^{\infty}_{0}(\bR^{d+1})$,
$r>0$, and $\nu\geq 16$, satisfying
$\nu r\leq1$, we have
\begin{equation}
                                                   \label{06.7.27.1}
\dashint_{(0,r^{2})}\dashint_{(0,r^{2})}|\,\|u_{xx}(t,\cdot)
\|_{L_{p}(\bR^{d})}-\|u_{xx}(s,\cdot)
\|_{L_{p}(\bR^{d})}|^{p}\,dtds
\end{equation}
$$
\leq N(\nu^{-p}
+\nu^{d+2}\hat{a}^{1/2})\dashint_{(0,\nu^{2}r^{2})}
\|u_{xx}(t,\cdot)
\|_{L_{p}(\bR^{d})}^{p}\,dt
$$
$$
+N\nu^{d+2}\dashint_{(0,\nu^{2}r^{2})}
\| Lu (t,\cdot)
\|_{L_{p}(\bR^{d})}^{p}\,dt.
$$

\end{corollary}

Indeed, by the triangle inequality
$$
|\|u_{xx}(t,\cdot)
\|_{L_{p}(\bR^{d})}-\|u_{xx}(s,\cdot)
\|_{L_{p}(\bR^{d})}|^{p}\leq
\|u_{xx}(t,\cdot)-u_{xx}(s,\cdot)\|_{L_{p}(\bR^{d})}^{p},
$$
so that the left-hand side of \eqref{06.7.27.1} is less than
$$
I:=\dashint_{(0,r^{2})}\dashint_{(0,r^{2})} 
\int_{\bR^{d}}|u_{xx}(t,x)-u_{xx}(s,x)|^{p}
\,dxdtds
$$
$$
=\dashint_{(0,r^{2})}\dashint_{(0,r^{2})} 
\int_{\bR^{d}}|u_{xx}(t,x+y)-u_{xx}(s,x+y)|^{p}
\,dxdtds,
$$
where $y$ is any point in $\bR^{d}$. By taking the average
of the extreme terms over $y\in B_{r}$ we see that
\begin{equation}
                                                \label{06.7.27.3}
I=\dashint_{(0,r^{2})}\dashint_{(0,r^{2})} 
\int_{\bR^{d}}\big(\dashint_{B_{r}(x)}
|u_{xx}(t,z)-u_{xx}(s,z)|^{p}\,dz\big)\,dxdtds.
\end{equation}

Next, since
$$
|u_{xx}(t,z)-u_{xx}(s,z)|^{p}\leq2^{p-1}
|u_{xx}(t,z)-(u_{xx})_{Q_{r}(0,x)} |^{p}
$$
$$
+2^{p-1}|u_{xx}(s,z)-(u_{xx})_{Q_{r}(0,x)}|^{p},
$$
we have that
$$
I\leq2^{p}\int_{\bR^{d}}(|u_{xx}-(u_{xx})_{Q_{r}(0,x)}|^{p}
)_{Q_{r}(0,x)}\,dx.
$$
By Lemma \ref{lemma 06.7.26.2} applied
to shifted cylinders the last expression is
dominated by a constant times
$$
(\nu^{-p}+\nu^{d+2}\hat{a}^{1/2})\int_{\bR^{d}}(|u_{xx} |^{p}
)_{Q_{\nu r}(0,x)}\,dx
+\nu^{d+2}\int_{\bR^{d}}(|Lu  |^{p}
)_{Q_{\nu r}(0,x)}\,dx,
$$
which similarly to \eqref{06.7.27.3} is shown to equal
$$
(\nu^{-p}+\nu^{d+2}\hat{a}^{1/2})\dashint_{(0,\nu^{2}r^{2})}
\|u_{xx}(t,\cdot)
\|_{L_{p}(\bR^{d})}^{p}\,dt
$$
$$+ \nu^{d+1}\dashint_{(0,\nu^{2}r^{2})}
\| Lu (t,\cdot)
\|_{L_{p}(\bR^{d})}^{p}\,dt
$$
and this yields \eqref{06.7.27.1}.

To move further  fix   
a    $u\in C^{\infty}_{0}(\bR^{d+1})$ and set
$$
\phi(t)=\|u_{xx}(t,\cdot)\|_{L_{p}(\bR^{d})},\quad
f=Lu ,\quad
\psi(t)=\|f(t,\cdot)\|_{L_{p}(\bR^{d})}
$$
and for any locally integrable function
$\tau(s)$ on $\bR$ denote by
$$
 M_{t}\tau(s)\quad\text{and}\quad \tau^{\#(t)}(s)
$$ 
the maximal and sharp functions of $\tau$, respectively.
\begin{lemma}
                                           \label{lemma 06.7.28.2}
Let $r_{0}\in(0,\infty)$,  
 $b=0$, $c=0$. Assume that the above $u(t,x)
=0$ for $t\not\in(0,r^{2}_{0})$. Then
for any $\nu\geq 16$ and
$R \in(0,1  ]$, we have  
$$
\phi^{\#(t)}\leq
N\nu^{(d+2)/p}M^{1/p}_{t}(\psi^{p})
$$
\begin{equation}
                                               \label{06.7.28.6}
+N\big((\nu r_{0}/R)^{2-2/p}
+\nu^{-1}
+\nu^{(d+2)/p}\omega^{1/(2p)}(R)\big)
M_{t}^{1/ p }(\phi^{p}),
\end{equation}
where $N=N(\omega,d,\kappa,K,p)$.
\end{lemma}

Proof. Obviously, Corollary \ref{corollary 06.7.28.1} in terms of
the functions $\phi$ and $\psi$ yields
$$
\dashint_{(0,r^{2})}\dashint_{(0,r^{2})}|\phi(t)-\phi(s)|^{p}
\,dtds
\leq
 N\nu^{d+2}\dashint_{(0,\nu^{2}r^{2})}\psi^{p}(t)\,dt
$$
$$
+ N(\nu^{-p}+\nu^{d+2}\omega ^{1/2}(R))
\dashint_{(0,\nu^{2}r^{2})}\phi^{p}(t)\,dt
$$
if $r\leq R/\nu$ (when $a^{\#(x)}_{\nu r}\leq
a^{\#(x)}_{R}\leq\omega(R)$ and $\nu r\leq1$). This corollary 
allows shifting  the origin.
Therefore, for any $\alpha,\beta\in\bR$ such that
$\alpha<\beta$ and $\beta-\alpha=r^{2}\leq R^{2}/\nu^{2}$ we have
$$
\dashint_{(\alpha,\beta)}\dashint_{(\alpha,\beta)}
|\phi(t)-\phi(s)|^{p}
\,dtds
\leq
 N\nu^{d+2}\dashint_{(\alpha,\alpha+
\nu^{2}(\beta-\alpha))}\psi^{p}(t)\,dt
$$
$$
+ N(\nu^{-p}+\nu^{d+2}\omega ^{1/2}(R))
\dashint_{(\alpha,\alpha+\nu^{2}(\beta-\alpha))}\phi^{p}(t)\,dt.
$$
Take a point $t_{0}\in\bR$ and $\alpha$ and $\beta$ as above
and such that $t_{0}\in(\alpha,\beta)$. Then
$t_{0}\in (\alpha,\alpha+\nu^{2}(\beta-\alpha))$ and by definition
$$
\dashint_{(\alpha,\alpha+
\nu^{2}(\beta-\alpha))}\psi^{p}(t)\,dt\leq M_{t}
(\psi^{p})(t_{0}),
$$
$$
\dashint_{(\alpha,\alpha+\nu^{2}(\beta-\alpha))}\phi^{p}(t)\,dt
\leq M_{t}(\phi^{p})(t_{0}) .
$$
By applying  H\"older's inequality   we conclude that
\begin{equation}
                                               \label{06.7.31.4}
\dashint_{(\alpha,\beta)}\dashint_{(\alpha,\beta)}|\phi(t)-\phi(s)| 
\,dtds
\end{equation}
is dominated by the value at $t_{0}$ of the right-hand
side of \eqref{06.7.28.6}, whenever 
$t_{0}\in(\alpha,\beta)$ and $\beta-\alpha\leq
R^{2}/\nu^{2}$. However,
if $\beta-\alpha>
R^{2}/\nu^{2}$, then \eqref{06.7.31.4} is dominated by
$$
2\dashint_{(\alpha,\beta)}  
I_{(0,r_{0}^{2})}\phi\,dt\leq
2\big(\dashint_{(\alpha,\beta)}  
I_{(0,r_{0}^{2})} \,dt\big)^{1-1/p}
\big(\dashint_{(\alpha,\beta)}  
 \phi^{p}\,dt\big)^{ 1/p}
$$
$$
\leq2(r_{0}^{2}/(\beta-\alpha))^{1-1/p}M^{1/p}_{t}
(\phi^{p})(t_{0})
\leq2(\nu r_{0} /R)^{2-2/p}M^{1/p}_{t}
(\phi^{p})(t_{0}).
$$
In this case \eqref{06.7.31.4} is again less than
the value at $t_{0}$ of
the right-hand
side of \eqref{06.7.28.6}. By taking the supremum
of  \eqref{06.7.31.4} over all $\alpha<\beta$ such that
$t_{0}\in(\alpha,\beta)$ we obtain \eqref{06.7.28.6}
at $t_{0}$. Since $t_{0}$ is arbitrary, the
lemma is proved.

\begin{lemma}
                                     \label{lemma 06.8.2.2}
There exists a constant 
  $N$ depending
only on $p,q,d,\kappa,K$, and the function
 $\omega$, such that   for any $u\in C^{\infty}_{0}(\bR^{d+1})$,
\begin{equation}
                                           \label{06.8.2.2}
 \|u_{xx}\|_{L_{q,p}}+\|u_{t}\|_{L_{q,p}}\leq
N(\|Lu \|_{L_{q,p}}+
 \|u_{x}\|_{L_{q,p}}+
 \|u\|_{L_{q,p}}).
\end{equation}

\end{lemma}

Proof. 
Notice that we included $\|u_{x}\|_{L_{q,p}}$
and $\|u \|_{L_{q,p}}$ on the right. Therefore,
while proving \eqref{06.8.2.2} we
  may certainly assume that $b\equiv0$ and $c\equiv0$.
Since $u_{t}= Lu 
-a^{ij}u_{x^{i}x^{j}}$, we only need to estimate $u_{xx}$.
If $p=q$ so that
$ L_{q,p} =L_{p} $, the result is known
from \cite{Kr06}. 

In case $q>p$ we fix a number $r_{0}$ and 
first assume that 
$$
u(t,x)=0\quad\text{for}\quad
 t\not\in(0,r^{2}_{0}).
$$
 Then set $f=Lu $
and also use other objects
introduced before Lemma \ref{lemma 06.7.28.2}.
We raise both parts of \eqref{06.7.28.6} to the power
$q$, integrate over $\bR$, and observe that since $q/p>1$,
by the Hardy-Littlewood theorem we have
$$
\int_{\bR}M_{t}^{q/p}(\psi^{p})(t)\,dt
\leq N\int_{\bR}\psi^{q}(t)\,dt=N\|f\|^{p}_{L_{q,p} },
$$ 
$$
\int_{\bR}M_{t}^{q/p}(\phi^{p})(t)\,dt
\leq  N\|u_{xx}\|^{q}_{L_{q,p}}.
$$ 
We also use the Fefferman-Stein theorem and conclude that
$$
 \|u_{xx}\|_{L_{q,p}} \leq
N_{1}\nu^{(d+2)/p} \|f\|_{L_{q,p} }
$$
\begin{equation}
                                           \label{06.8.2.3} 
+N_{2}\big((\nu r_{0}/R)^{2-2/p}
+\nu^{-1}
+\nu^{(d+2)/p}\omega^{1/(2 p)}(R)\big)
\|u_{xx}\|_{L_{q,p} }),
\end{equation}
whenever $\nu\geq16$ and $R\leq1$, where
$N_{i}$ are determined by $p,q,d,\kappa,K$ and the  
function $\omega$. We choose a large $\nu=\nu(N_{2},d)$ 
and a small $R=R(N_{2},d,q,\omega)$ so that
$$
N_{2}\big( \nu^{-1}
+\nu^{(d+2)/p}\omega^{1/(2 p)}(R)\big)\leq1/4.
$$
After $\nu$ and $R$ have been fixed, we chose a small $r_{0}
=r_{0}(N_{2},d,q,\omega)$
so that
$$
N_{2} (\nu r_{0}/R)^{2-2/q}\leq1/4.
$$
Then \eqref{06.8.2.3}  implies that
\begin{equation}
                                           \label{06.8.2.4} 
\|u_{xx}\|_{L_{q,p} } \leq2
N_{1}\nu^{(d+2)/p} \|Lu \|_{L_{q,p} }
\end{equation}
for any $u\in C^{\infty}_{0}(\bR^{d+1})$
such that $u(t,x)=0$ if $t\not\in(0,r_{0}^{2})$.
We thus have obtained \eqref{06.8.2.2}
  even without the terms
$\|u_{x}\|_{L_{q,p} }$
and $\|u \|_{L_{q,p} }$ on the right of~\eqref{06.8.2.2}.

Now take a nonnegative $\zeta\in C^{\infty}_{0}(\bR)$ such that
$\zeta(t)=0$ if $t\not\in(0,r_{0}^{2})$ and
$$
\int_{\bR}\zeta^{p}(t)\,dt=1.
$$
Also take a $u\in C^{\infty}_{0}(\bR^{d+1})$ and observe that
\eqref{06.8.2.4}  is also true if we shift the $t$ axis.
In particular, \eqref{06.8.2.4}  is applicable
to $u(t,x)\zeta(t-t_{0})$.
Then we get
$$
\int_{\bR}\zeta^{q}(t-t_{0})
\|u_{xx}(t,\cdot)\|^{q}_{L_{p}(\bR^{d})}\,dt
\leq N 
\int_{\bR}\zeta^{q}(t-t_{0})
\| Lu (t,\cdot)\|^{q}_{L_{p}(\bR^{d})}\,dt
$$
$$
+N\int_{\bR}
|\zeta'(t-t_{0})|^{q}
\|u (t,\cdot)\|^{q}_{L_{p}(\bR^{d})}\,dt 
$$
Upon integrating through with respect to $t_{0}$
we come to \eqref{06.8.2.2}. The lemma is proved.

On the basis of this lemma by repeating
almost word for word the proof
of Theorem 4.1 of \cite{Kr06}
(or using the method of proving Theorem \ref{theorem 7.19.1}
or Lemma  \ref{lemma 06.5.18.1})
we obtain the following result.
\begin{theorem}
                                             \label{theorem 06.10.17.5}
There are   constants
$\lambda_{0}$ and
$N$, depending only on $p$, $K$, $\kappa$, $d$, and 
$\omega$, such that for any $\lambda\geq\lambda_{0}$
and $u\in W^{1,2}_{q,p}$ we have  
$$
\lambda\|u\|_{L_{q,p}}+\sqrt{\lambda}\|u_{x }\|_{L_{q,p}}
+\|u_{xx}\|_{L_{q,p}}
+\|u_{t}\|_{L_{q,p}}\leq N\|(L-\lambda)u\|_{L_{q,p}}.
$$
 
Furthermore, for any $\lambda\geq\lambda_{0}$ and
$f\in L_{q,p}$ there exists a unique 
$u\in W^{1,2}_{q,p}$ such that $(L-\lambda)u=f$.

\end{theorem}

Finally, Theorem \ref{theorem 06.10.17.5} implies
Theorem \ref{theorem 7.20.3} in the same way as Theorem
4.1 of \cite{Kr06} implies Theorem 2.1 of \cite{Kr06}.

\mysection{Proof of Theorem \protect\ref{theorem 7.20.4}}
                                     \label{section 06.10.22.2}

We start with the following result which 
will be proved in Section \ref{section 06.10.20.1}
and which
is
an improvement of the key estimate found in the proof
of Theorem 5.3 of \cite{Kr06}.
We work in the setting in which Theorem \ref{theorem 7.20.4}
is stated.

\begin{lemma}
                                             \label{lemma 06.10.9.5}
Let  $b=\hat{b}=0$, $c=0$,
$f=(f^{1},...,f^{d})\in L_{p,loc} $.
Then 
  there exists a constant $N=N(d,\kappa,p,K, \omega)$ such that
for any $u\in \cH^{1}_{p,loc}$, $\nu\geq16$, and
$r\in(0,1/\nu]$, such that $\cL u=\Div f$ in $Q_{\nu r}$,
     we have
\begin{equation}
                                                \label{06.10.9.7}
 (|u_{x }-(u_{x })_{Q_{r} }|^{p})_{Q_{r} }
\leq N \nu ^{d+2} \cA_{\nu r}+
  N(\nu^{-p}+\nu ^{d+2} 
  \hat{a}^{1/2})\cB_{\nu r},
\end{equation}
where
$$
\cA_{\rho}= (
|f|^{p})_{Q_{\rho} },
\quad
\cB_{\rho}=( |u_{x }|^{p} )_{Q_{\rho} },
\quad
\hat{a}= a^{\#(x)}_{  \nu r }  .
$$
 \end{lemma}

The following is proved in the same way as
Corollary \ref{corollary 06.7.28.1}.

\begin{corollary}
                                         \label{corollary 06.10.19.1}
Let  $b=\hat{b}=0$, $c=0$,
$u\in \cH^{1}_{p }((S,T))$
for any finite $S<T$, $\cL u=\Div f$, where
$f=(f^{1},...,f^{d})\in L_{p}((S,T)\times\bR^{d}) $
for any finite $S<T$.
Then 
  there exists a constant $N=N(d,\kappa,p,K, \omega)$ such that
for any   $\nu\geq16$  and
$r\in(0,1/\nu]$ 
     we have
\begin{equation}
                                                   \label{06.10.19.1}
\dashint_{(0,r^{2})}\dashint_{(0,r^{2})}|\,\|u_{ x}(t,\cdot)
\|_{L_{p}(\bR^{d})}-\|u_{ x}(s,\cdot)
\|_{L_{p}(\bR^{d})}|^{p}\,dtds
\end{equation}
$$
\leq N(\nu^{-p}
+\nu^{d+2}\hat{a}^{1/2})\dashint_{(0,\nu^{2}r^{2})}
\|u_{ x}(t,\cdot)
\|_{L_{p}(\bR^{d})}^{p}\,dt
$$
$$
+N\nu^{d+2}\dashint_{(0,\nu^{2}r^{2})}
\| f(t,\cdot)
\|_{L_{p}(\bR^{d})}^{p}\,dt.
$$

\end{corollary}

After that in the same way as Lemma \ref{lemma 06.8.2.2}
is proved one derives its counterpart
for divergence equations from 
Corollary \ref{corollary 06.10.19.1}.

\begin{lemma}
                                        \label{lemma 06.10.19.1}

Let $q\geq p$, $u\in\cH^{1}_{q,p}$,
$b=\hat{b}=0$, $c=0$, $\cL u=\Div f$ with
$f\in L_{q,p}$. Then there exists a constant $N$,
depending only on $q,p,d,\kappa,K$, and $\omega$ such that 
\begin{equation}
                                                   \label{06.10.19.2}
\|u_{x}\|_{L_{q,p}}\leq N(\|f\|_{L_{q,p}}
+\|u\|_{L_{q,p}}).
\end{equation}
\end{lemma}

Next we state and  prove
an analog of Lemma 5.5  of \cite{Kr06} where $q=p$
and $u$ is supposed to have small support.

\begin{theorem}
                                             \label{theorem 7.19.1}
Let $q\geq p$, 
$f=(f^{1},...,f^{d})$,
$f^{i}, g\in L_{q,p} $,
$u\in \cH^{1}_{q,p} $, $\lambda\in\bR$, and
$$
\cL u-\lambda u=\Div f+g.
$$
 We assert that
 there exist  constants 
 $\lambda_{0},N\in(0,\infty)$,
 depending only on
 $p,q$, $d$, $K$, $\kappa$, and $\omega$,
such that
\begin{equation}
                                                \label{7.19.2}
\|u_{t}\|_{\bH^{-1}_{q,p}}+
\sqrt{\lambda}\|u_{x }\|_{L_{q,p}}+
\lambda \|u \|_{L_{q,p}}\leq N(\sqrt{\lambda}\|f\|_{L_{q,p}}
+\|g\|_{L_{q,p}}),
\end{equation}
provided that   $\lambda\geq\lambda_{0}$.
\end{theorem}

Proof. 
We follow the same pattern as in the proof of Lemma 5.5 of \cite{Kr06}.
 First, we observe that
the terms $(\hat{b}^{i}u)_{x^{i}}$ and
$b^{i}u_{x^{i}}+cu$ in $\cL u$ can be included in $\Div f$
and $g$, respectively. This will introduce new terms
in the right-hand side of \eqref{7.19.2}
but on the account
of perhaps increasing $\lambda_{0}$ they can be
absorbed into the left-hand side of \eqref{7.19.2}.
For this reason in the rest of the proof
we may and will assume that $b=\hat{b}=0$, $c=0$.

In this case we use a method introduced by Agmon.
Consider the space
$\bR^{d+2}=\{(t,z)=(t,x,y):t,y\in\bR,x\in\bR^{d}\}$ and the function
\begin{equation}
                                                     \label{7.19.3}
\tilde{u}(t,z)=u(t,x)\zeta(y)\cos(\mu y),
\end{equation}
 where
 $\mu=\sqrt{\lambda}$ and $\zeta$ is an
odd $C^{\infty}_{0}(\bR)$-function,
$\zeta\not\equiv0$.
Also introduce the operator
\begin{equation}
                                                     \label{7.19.4}
\tilde{\cL}u(t,z)=\cL(t,x)u(t,z)+u_{yy}(t,z).
\end{equation}
As in \cite{Kr06} one checks that the coefficients
of $\tilde{\cL}$ are $VMO_{x}$-functions
(with respect to $(t,z)$).

Set  
$\tilde{f}^{i}(t,z)=f^{i}(t,x)\zeta(y)\cos(\mu y)$ for
$i=1,...,d$ and
$$
\tilde{f}^{d+1}(t,z)=g(t,x)\zeta_{1}(y)
-2u(t,x)\zeta_{2}(y) +
u(t,x)\zeta_{3}(y),
$$
where
$$
\zeta_{1}(y)=\int_{-\infty}^{y}\zeta(s)\cos(\mu s)\,ds,\quad
\zeta_{3}(y)=\int_{-\infty}^{y}\zeta''(s)\cos(\mu s)\,ds
$$
$$
\zeta_{2}(y)=\mu\int_{-\infty}^{y}\zeta'(s)\sin(\mu s)\,ds
=-\zeta'(y)\cos(\mu y)+\zeta_{3}(y).
$$
Observe that $\zeta_{i}\in C^{\infty}_{0}(\bR)$ since
$\zeta$ is odd and has compact support. Furthermore, as is easy to
check,
$$
\tilde{\cL} \tilde{u}(t,z)
=(\tilde{f}^{1}(t,z))_{x^{1}}+...+
(\tilde{f}^{d}(t,z))_{x^{d}}+(\tilde{f}^{d+1}(t,z))_{y}.
$$

We denote by $\tilde{L}_{p}$ the $L_{p}$ space of functions
of $z=(x,y)$ (avoiding using a confusing notation $L_{p}(\bR^{d+1})$)
and by Lemma \ref{lemma 06.10.19.1} obtain
\begin{equation}
                                                \label{7.19.5}
\int_{\bR}\|\tilde{u}_{z}(t,\cdot)\|^{q}_{\tilde{L}_{p}}\,dt
\leq
N\sum_{i=1}^{d+1}\int_{\bR}
\|\tilde{f}^{i}(t,\cdot)\|^{q}_{\tilde{L}_{p}}\,dt
+N\int_{\bR}\|\tilde{u}(t,\cdot)\|^{q}_{\tilde{L}_{p}}\,dt.
\end{equation}

Since
$$
\delta_{0}:=\int_{\bR^{d}}|\zeta(y)\sin(\mu y)|^{p}\,dy,\quad
\delta_{1}:=\int_{\bR^{d}}|\zeta(y)\cos(\mu y)|^{p}\,dy
$$
are bounded  away from zero for   $\mu\geq1$, we get
for each $t$ and $\mu\geq1$ that
$$
\|u_{x }(t,\cdot)\|^{p}_{L_{p}(\bR^{d}) }
=\delta_{1}^{-1}\int_{\bR^{d+1}}|u_{x}(t,x)
\zeta (y)\cos(\mu y)|^{p}
\,dz \leq\delta_{1}^{-1}\|\tilde{u}_{z}(t,\cdot)\|^{p}_{\tilde{L}_{p}},
$$ 
$$
\|u (t,\cdot)\|^{p}_{L_{p}(\bR^{d}) }=
\delta_{0}^{-1}\mu^{-p}\int_{\bR^{d+1}}|\tilde{u}_{y}(t,z)-
u(t,x) \zeta'(y)\cos(\mu y)|^{p}\,dz 
$$
$$
\leq N\mu^{-p}(\|\tilde{u}_{z}(t,\cdot)\|^{p}_{\tilde{L}_{p}}
+\|u(t,\cdot) \|^{p}_{L_{p}(\bR^{d}) }).
$$
It follows that if $\mu$ is large enough, then
$$
\mu^{p}\|u (t,\cdot)\|^{p}_{L_{p}(\bR^{d}) }
\leq N\|\tilde{u}_{z}(t,\cdot)\|^{p}_{\tilde{L}_{p}}.
$$

Hence, by \eqref{7.19.5}   for large $\mu$  
\begin{equation}
                                                \label{7.19.6}
\mu^{q}\|u \|^{q}_{L_{q,p} }+\|u_{x }\|^{q}_{L_{q,p} }
\leq N\sum_{i=1}^{d+1}\int_{\bR}
\|\tilde{f}^{i}(t,\cdot)\|^{q}_{\tilde{L}_{p}}\,dt
+N\int_{\bR}\|\tilde{u}(t,\cdot)\|^{q}_{\tilde{L}_{p}}\,dt.
\end{equation}

Now we estimate the right-hand side of \eqref{7.19.6}.
Obviously, 
$$
\|\tilde{f}^{i}(t,\cdot)\|^{q}_{\tilde{L}_{p}}\leq N
\|f^{i}(t,\cdot) \|^{q}_{ L_{p}(\bR^{d}) }, \quad
i=1,...,d,
$$
$$
\|\tilde{u} (t,\cdot)\|^{q}_{\tilde{L}_{p}}\leq N
\|u (t,\cdot) \|^{q}_{ L_{p}(\bR^{d}) }.
$$

Furthermore,  
$$
\zeta_{1}=\mu^{-1}\big[\zeta (y)\sin(\mu y)-\int_{-\infty}^{y}
\zeta' (s)\sin(\mu s)\,ds\big],
$$
which shows that $\zeta_{1}$ equals $\mu^{-1}$
times a uniformly bounded function with support not wider than that
of $\zeta$. Hence,
$$
\|g\zeta_{1}(t,\cdot)\|^{q}_{\tilde{L}_{p}}\leq N\mu^{-q}
\|g(t,\cdot) \|^{q}_{ L_{p}(\bR^{d}) }.
$$
Also $\zeta_{2}$ and $\zeta_{3}$ are  uniformly bounded 
 with support not wider than that
of $\zeta$. Therefore,
$$
\|(2u\zeta_{2}-u\zeta_{3})(t,\cdot)\|^{q}_{\tilde{L}_{p}}\leq N
\|u (t,\cdot) \|^{q}_{ L_{p}(\bR^{d}) },
$$
$$
\|\tilde{f}^{d+1} (t,\cdot)\|^{q}_{\tilde{L}_{p}}\leq N\mu^{-q}
\|g(t,\cdot) \|^{q}_{ L_{p}(\bR^{d}) }
+N\|u (t,\cdot) \|^{q}_{ L_{p}(\bR^{d}) }.
$$
This and \eqref{7.19.6} yield \eqref{7.19.2} without the term with
$u_{t}$. To estimate this term it suffices to observe that
$$
(1-\Delta)^{-1/2}u_{t}=-(1-\Delta)^{-1/2}D_{j}(a^{ij}u_{x^{i}}-f^{j})
+(1-\Delta)^{-1/2}(\lambda u+g),
$$
so that, by the boundedness of $(1-\Delta)^{-1/2} $
and $(1-\Delta)^{-1/2}D_{j}$, for each~$t$
$$
\|(1-\Delta)^{-1/2}u_{t}(t,\cdot)\|_{L_{p}(\bR^{d})}
\leq N(\|u_{x}(t,\cdot)\|_{L_{p}(\bR^{d})}
+\lambda\|u (t,\cdot)\|_{L_{p}(\bR^{d})}
$$
$$
+\|f(t,\cdot)\|_{L_{p}(\bR^{d})}+\|g(t,\cdot)\|_{L_{p}(\bR^{d})})
$$
Upon raising both parts to the power $q$ and integrating
over $t\in\bR$ we get the required estimate of $u_{t}$.
The theorem is proved.

A simple argument in Section 6 of \cite{Kr06}
shows that
$$
\|u\|_{L_{q,p}}+\|u_{x}\|_{L_{q,p}}+\|u_{t}\|_{\bH^{-1}_{q,p}} 
\quad\text{and}\quad
\|(1-\Delta)^{-1/2}u\|_{W^{1,2}_{q,p}}  
$$
define   equivalent norms in $\cH^{1}_{q,p}$.
This argument also shows that, for each fixed $\lambda>0$,
the right-hand side of \eqref{7.19.2} dominates
$$
\|\Div f+g\|_{\bH^{-1}_{q,p}}
$$
and in turn one can find $\tilde{f}$ and $\tilde{g}$
so that $\Div f+g=\Div\tilde{f}+\tilde{g}$ and the 
right-hand side of \eqref{7.19.2} is dominated by
$$
N\|\Div\tilde{f}+\tilde{g}\|_{\bH^{-1}_{q,p}}.
$$

Therefore, Theorem \ref{theorem 7.19.1} implies
assertion (i) for $q\geq p$ in the following result.

\begin{theorem}
                                       \label{theorem 7.22.3}
  There is a constant
$\lambda_{0}$ depending only on $p,q$, $d$, $\kappa$,
$K$, and $\omega$ such that for any $\lambda\geq\lambda_{0}$

(i) for
  any $u\in\cH^{1}_{q,p}$ we have
\begin{equation}
                                                \label{06.10.19.7}
\|u\|_{\cH^{1}_{q,p}}\leq N(\lambda,p,d,\kappa,K,\omega)
\|(\cL-\lambda)u\|_{\bH^{-1}_{q,p}};
\end{equation}

(ii) for any $h\in\bH^{-1}_{q,p}$ there exists a unique
$u\in\cH^{1}_{q,p}$ such that $\cL u-\lambda u=h$.
\end{theorem}

Proof. 
It is a classical  result that for any $\lambda>0$
and $g\in L_{q,p}$ there exists a (unique) solution
$w\in W^{1,2}_{q,p}$
of $\Delta w+w_{t}-\lambda w=g$ and one even can give
$w$ by a formula (see, for instance,
  Theorem 4.2 of \cite{Kr02} and the references
in \cite{Kr02}). Then
$u:=(1-\Delta)^{1/2}w$ is in $\cH^{1}_{q,p}$ and satisfies
$\Delta u+u_{t}-\lambda u=h$
with $h=(1-\Delta)^{1/2}g$.
As $g$ runs through $L_{q,p}$, $h$ runs through $\bH^{-1}_{q,p}$
by definition.

Hence, the present theorem holds if $\cL u=\Delta u+u_{t}$. By what
has been said before the theorem the a priori estimate 
\eqref{06.10.19.7} holds if $q\geq p$. Then by the method of
continuity assertion (ii) also holds if $q\geq p$.

The case $1<q<p$ is considered in a standard way by duality
owing to the fact that the formally adjoint operator
to $\cL$ has the same structure as $\cL$ 
only with reversed time
axis.   
The theorem is proved.

Finally, Theorem \ref{theorem 7.20.3} is derived from Theorem
\ref{theorem 7.22.3}
  in the same way as in similar situations in \cite{Kr06}.

\mysection{New approach to the $L_{p}$-theory for equations  
with VMO coefficients}  
                                       \label{section 7.22.2}

We
take an $a\in\bA$ and set 
$$
\bar{L}u(t,x)=a^{ij}(t)u_{x^{i}x^{j}}(t,x)
+u_{t}(t,x).
$$ 
In this section $p\in(1,\infty)$ and $\lambda\geq0$
unless explicitly specified otherwise.

Here we give a new proof of the following result from
\cite{Kr06}, which is a simplified version of
Lemma \ref{lemma 06.7.26.2}.
\begin{theorem}
                           \label{theorem 06.6.13.1} 
 
There is a constant
$N$, depending only on $d,p,K$, and $\kappa$, such that
for any $u\in W^{1,2}_{p,loc} $,
$r\in(0,\infty)$, and $\nu\geq4$
\begin{equation}
                                                     \label{06.6.13.10}
 (|u_{xx}(t,x) 
-(u_{xx})_{Q_{r}}|^{p})_{Q_{r}}\leq
N\nu^{d+2}(|\bar{L}u |^{p})_{Q_{\nu r}}
+ N\nu^{-p}
(|u_{xx}|^{p})_{Q_{\nu r}}.
\end{equation}
\end{theorem}

In \cite{Kr06} Theorem \ref{theorem 06.6.13.1}
is proved
on the basis of solving boundary-value problems
for parabolic equations. The proof we give later in the section
is based on solvability of equations in the whole space
and
extends to more general operators and systems of equations
without much effort. 
In particular, we will see that,
once the solvability theory for
operators $\bar{L}$ is developed in $W^{1,2}_{p}(\bR^{d+1})$
for a $p>1$, Theorem \ref{theorem 06.6.13.1}
becomes available and, according to simple arguments
from \cite{Kr06}, the solvability theory
in $W^{1,2}_{q}$ with $q>p$ for equations with $VMO_{x}$ coefficients
becomes available as well.

This fact has the following methodological
implication.
If $p=2$ one can construct 
the solvability theory for $\bar{L}$ in $W^{1,2}_{2}(\bR^{d+1})$
by using the Fourier transform. Then by the above
(or by what is done in Remark \ref{remark 06.10.17.3} below),
the solvability theory for $\bar{L}$ in $W^{1,2}_{p}(\bR^{d+1})$
with $p>2$ is available. By duality one gets it for
$p\in(1,2)$ as well and as has been pointed out,
this is the only thing one needs to 
construct the solvability theory for
operators with $VMO_{x}$ coefficients in $W^{1,2}_{p}(\bR^{d+1})$,
$p\in(1,\infty)$.
 
As usual, for any multi-index $\alpha=(\alpha_{1},...,\alpha_{d})$,
$\alpha_{i}\in\{0,1,2,...\}$, we set
$$
D^{\alpha}u=D^{\alpha_{1}}_{1}\cdot...\cdot
D^{\alpha_{d}}_{d}u,\quad D_{i}u=u_{x^{i}}=\frac{\partial u}
{\partial x^{i}},\quad|\alpha|=\alpha_{1}+...+\alpha_{d}.
$$

\begin{lemma}
                                                \label{lemma 9.15.1}
Take   $p\in[1,\infty)$  and
$N_{0}\in(0,\infty)$  and
assume that  for any $u\in W^{1,2}_{p}(\bR^{d+1}_{0})$
 we have
\begin{equation}
                                                  \label{9.15.2} 
 \|u_{t}\|_{L_{p}(\bR^{d+1}_{0})} +
\|u_{xx}\|_{L_{p}(\bR^{d+1}_{0})}
\leq N_{0}(\|Lu \|_{L_{p}(\bR^{d+1}_{0})}+\| u
\|_{L_{p}(\bR^{d+1}_{0})}).
\end{equation}

Then  for any $0<r<R<\infty$
there exists a constant $N$, depending only
on $N_{0},d,p,K,r$, and $ R  $,
such that for any   $u\in
W^{1,2}_{p}(Q_{R})$ we have
\begin{equation}
                                                  \label{9.15.3}
\|u_{t}\|_{L_{p}(Q_{r})}
+\|u_{xx}\|_{L_{p}(Q_{r})}\leq N\big(\|
Lu \|_{ L_{p}(Q_{R})}
+ \|u_{x}\|_{L_{p}(Q_{R})}+ \|u\|_{L_{p}(Q_{R})}\big).
\end{equation}
\end{lemma}

This is a trivial result, which is obtained by taking an appropriate
cut-off function $\zeta$ and applying \eqref{9.15.2}
to $u\zeta$.
\begin{remark}
With a little extra work
(see the proof of Lemma 4.2 of \cite{KDK2})
 one shows that the term with $u_{x}$
on the right in \eqref{9.15.3} can be dropped.
\end{remark}

Another general result we need  is
a parabolic analog of Poincar\'e's
inequality. It is proved in the same way as Lemma 3.2 of \cite{Kr06}
(also see Lemma 4.2 of \cite{PS}).
We generalize Lemma \ref{lemma 06.7.7.2} in Lemmas
\ref{lemma 06.7.17.1} and \ref{lemma 06.7.26.1}.

\begin{lemma}
                                            \label{lemma 06.7.7.2}
Let $p\in[1,\infty)$. Then 
there is a constant $N=N(d,p )$ such that
for any $r\in(0,\infty)$ and 
  $u\in
C^{\infty}_{loc}(\bR^{d+1})$ we have
\begin{equation}
                                                     \label{06.7.7.6}
\int_{Q_{r}}|u_{x }(t,x) 
-(u_{x })_{Q_{r}}|^{p} \,dxdt
\leq Nr^{p}\int_{Q_{r}}(|u_{xx}|^{p} +|u_{t}|^{p} )\,dxdt,
\end{equation}
$$
\int_{Q_{r}}|u(t,x)-u_{Q_{r}}-
x^{i}(u_{x^{i}})_{Q_{r}}|^{p} \,dxdt
$$
\begin{equation}
                                                     \label{06.7.7.2}
\leq Nr^{2p}\int_{Q_{r}}(|u_{xx}|^{p} +|u_{t}|^{p} )\,dxdt.
\end{equation}

 \end{lemma}

We need the following classical result
(which can be obtained, for instance, along the lines discussed
 after Theorem \ref{theorem 06.6.13.1}).
\begin{theorem}
                                            \label{theorem 06.10.6.2}
There is a constant $N=N(p,d,\kappa,K)$
such that for any $\lambda\geq 0$, $T\in[-\infty,\infty)$,
and $u\in W^{1,2}_{p}(\bR^{d+1}_{T})$ we have
$$
\lambda \|u\|_{L_{p}(\bR^{d+1}_{T})}
+\|u_{xx}\|_{L_{p}(\bR^{d+1}_{T})}+\|u_{t}\|_{L_{p}(\bR^{d+1}_{T})}
\leq N\|\bar{L}u-\lambda u\|_{L_{p}(\bR^{d+1}_{T})}.
$$
Furthermore, for any $\lambda>0$ and $f\in L_{p}(\bR^{d+1}_{T})$
there exists a unique $u\in W^{1,2}_{p}(\bR^{d+1}_{T})$
such that $\bar{L}u+u_{t}-\lambda u=f$.
\end{theorem}

\begin{remark}
                                       \label{remark 06.10.17.1}

Owing to Theorem \ref{theorem 06.10.6.2},
the assertion of
Lemma \ref{lemma 9.15.1} 
holds with $\bar{L}$ in place of $L$.
\end{remark}

\begin{remark}
                                       \label{remark 06.10.21.1}
In the proof of Theorem \ref{theorem 06.6.13.1}
we will use the decomposition $f:=\bar{L}u-\lambda u=g+h$,
where roughly speaking $g=fI_{Q_{\nu r}}$, and accordingly
have $u=v+w$, where $v$ is defined by the equation
$\bar{L}v-\lambda v=h$. The function $v$ is 
``harmonic" in $Q_{\nu r }$ in the sense that
$h=0$ there.
Then the oscillation of $w$ will be estimated by using 
Theorem \ref{theorem 06.10.6.2} and that of $v$
will be derived from what follows. Observe that
since we solve the equation $\bar{L}v-\lambda v=h$
in the whole space $\bR^{d+1}$ we need $\lambda>0$.
\end{remark}

\begin{lemma}  
                                              \label{lemma 06.5.27.1}
 Take $0<r<R<\infty$  and let
 $ m\in\{0,1,2,...\}$.
Take a function  $u\in W^{1, 2}_{p}(Q_{R})$
and assume that $ \bar{L}u  $ vanishes in $Q_{R}$.
 Then for any multi-index $\alpha$ the  derivatives
$D^{\alpha}u $ and $D^{\alpha}u_{t} $ are bounded in $Q_{r}$
and, with $N=N(|\alpha|,d,\kappa,K, r,R,p)$,   
\begin{equation}
                                                 \label{06.5.27.3}
\sup_{Q_{r} }
|D^{\alpha}u|
\leq N (\|u_{x}\|_{L_{p}(Q_{R})}+
\|u\|_{L_{p}(Q_{R})} )=:NI,\quad
 \sup_{Q_{r}  }
|D^{\alpha}u_{t}|
\leq NI .
\end{equation}
\end{lemma}

Proof. Since the coefficients of $\bar{L}$ are independent
of $x$ we can mollify the 
function $u$ with respect to $x$ and have
equation $\bar{L}\bar{u}=0$ 
in slightly smaller domain than $Q_{R}$
for $\bar{u}$ being the mollified $u$. Then, if
the result is true for $\bar{u}$, we can pass to the
limit as the support of the mollification kernel
shrinks to the origin. It follows that without losing generality
we may assume that $D^{\beta}u\in W^{1,2}_{p}(Q_{R})$
for any $\beta$. Then, since
$D^{\beta}u_{t}=-a^{ij}D^{\beta}u_{x^{i}x^{j}}$
in $Q_{R}$, we also have
 $D^{\beta}u_{t}\in W^{1,2}_{p}(Q_{R})$
for any $\beta$.
 
By Remark \ref{remark 06.10.17.1}, applied to $D^{\beta}u$,
for each integer
$k\geq0$ and $r<r_{1}<r_{2}<R$ we have  
$$
\sum_{|\beta|\leq k }\|D^{\beta}u_{t}\|_{L_{p}(Q_{r_{1}})}+
\sum_{|\beta|\leq k+2}\|D^{\beta}u\|_{L_{p}(Q_{r_{1}})}
\leq N
\sum_{|\beta|\leq k+1}\|D^{\beta}u\|_{L_{p}(Q_{r_{2}})},
$$
$$
\sum_{|\beta|\leq k+1}\|D^{\beta}u\|_{L_{p}(Q_{r_{1}})}
\leq N\bigg(
\sum_{|\beta|\leq k }\|D^{\beta}u\|_{L_{p}(Q_{r_{2}})}
+\|u_{x}\|_{L_{p}(Q_{r_{2}})}\bigg).
$$
By iterating the last relation we see that
$$
\sum_{|\beta|\leq k+1}\|D^{\beta}u\|_{L_{p}(Q_{s_{1} })}
\leq N(
 \| u\|_{L_{p}(Q_{s_{2}})}
+\|u_{x}\|_{L_{p}(Q_{s_{2}})})\leq NI,
$$
whenever $r<s_{1}<s_{2}<R$. Hence,
$$
\sum_{|\beta|\leq k }\|D^{\beta}u_{t}\|_{L_{p}(Q_{r })}+
\sum_{|\beta|\leq k+2}\|D^{\beta}u\|_{L_{p}(Q_{r })}
\leq NI.
$$
Furthermore, obviously
$$
|D_{t}\|D^{\beta}u(t,\cdot)\|_{L_{p}(B_{r})}|
\leq\|D^{\beta}u_{t}(t,\cdot)\|_{L_{p}(B_{r})}.
$$
Therefore,
for $\phi^{\beta}(t):=\|D^{\beta}u(t,\cdot)\|_{L_{p}(B_{r})}$
by embedding theorems we have
$$
\sup_{[0,r^{2}]}\phi^{\beta}\leq N(
\|\phi^{\beta}\|_{L_{p}(0,r^{2})}
+\|\phi^{\beta}_{t}\|_{L_{p}(0,r^{2})})
$$
$$
= N(\|D^{\beta}u\|_{L_{p}(Q_{r})}
+\|D^{\beta}u_{t}\|_{L_{p}(Q_{r})})\leq NI.
$$
Thus,
$$
\sup_{[0,r^{2}]} \sum_{|\beta|\leq k}
\|D^{\beta}u(t,\cdot)\|_{L_{p}(B_{r})}\leq NI.
$$
By embedding theorems, if $k$ is large enough, then
$$
\sup_{x\in B_{r}}|D^{\alpha}u(t,x)|
\leq N\sum_{|\beta|\leq k}
\|D^{\beta}u(t,\cdot)\|_{L_{p}(B_{r})}
$$
and this leads to the first estimate in \eqref{06.5.27.3}.
One gets the second one from the equation
$D^{\alpha}u_{t}=-a^{ij}D^{\alpha}u_{x^{i}x^{j}}$.
The lemma is proved. 

Below, for an integer $m\geq0$, by $D^{m}u(t,x)$
we mean the collection of all $m$th order derivatives
of $u$ with respect  to $x$. In the set
of these collection we define  a Euclidean norm
$|D^{m}u(t,x)|$.

\begin{lemma}
                                          \label{lemma 06.5.18.1}
Let   $m\in\{0,1,2,...\}$, $\lambda\geq0$,  and
 $u\in C^{\infty}_{0}(\bR^{d+1})$. Assume that 
$\bar{L}u -\lambda u$ vanishes in $Q_{2}$. Then,
with $N=N(d,m,\kappa,p,K)$,
\begin{equation}
                                                 \label{06.5.18.1}
\max_{Q_{1}}\big(|D^{m}u_{xx}|^{p}
+|D^{m}u_{t}|^{p}\big)
\leq N\int_{Q_{2}}(|u_{xx}|^{p} +|u_{t}|^{p}
+\lambda^{p/2}|u_{x}|^{p})\,dxdt.
\end{equation}
 
\end{lemma}

Proof. By  Lemma \ref{lemma 06.5.27.1} 
$$
I:=\max_{Q_{1}}\big(|D^{m}u_{xx}|^{p}
+|D^{m}u_{t}|^{p}\big)\leq N
(\|u_{x}\|^{p}_{L_{p}(Q_{3/2})}+\|u\|^{p}_{L_{p}(Q_{3/2})}).
$$
If  $\lambda=0$ we can replace here $u$ with $v:=u-u_{Q_{2}}
-x^{i}(u_{x^{i}})
_{Q_{2}}$ without violating the fact that $\bar{L}u+u_{t}$ 
vanishes in $Q_{2}$ or changing
the left-hand side. Therefore,  
$$
I\leq N(\|v_{x}\|^{p}_{L_{p}(Q_{2})}+\|v\|^{p}_{L_{p}(Q_{2})}),
$$
and using Lemma \ref{lemma 06.7.7.2}
  yields the desired result.

In the general case that $\lambda\geq0$ we again use
 a method 
suggested by S.~Agmon.
Introduce the function $\hat{u}(t,z)=\hat{u}(t,x,y)$ by
$$
\hat{u}(t,z)=u(t,x)\cos(\sqrt{\lambda}y)
$$
and set 
$$
\hat{Q}_{r}=(0,r^{2})\times\{|z|<r\}.
$$

Obviously,
$$
  D^{m}u_{xx}(t,x)= D^{m}\hat{u}_{xx}(t,x,0),\quad
 D^{m}u_{t}(t,x)= D^{m}\hat{u}_{t}(t,x,0)
$$
Therefore,
$$
I\leq
\max_{\hat{Q}_{1}}\big(|D^{m}\hat{u}_{xx}|^{p}+|D^{m}\hat{u}_{t}|^{p}\big).
$$
However,
$$
\bar{L}\hat{u}+\hat{u}_{yy} =0\quad\text{in}\quad\hat{Q}_{2},
$$
so that we can apply the above result to $\hat{u}$ and conclude
\begin{equation}
                                                 \label{06.5.18.5}
I
\leq N\int_{\hat{Q}_{2}}
\big(|\hat{u}_{zz}|^{p}+|\hat{u}_{t}|^{p}\big)\,dzdt.
\end{equation}
Here the term $\hat{u}_{zz}$ is the collection
consisting of 
$$
u_{xx}\cos(\sqrt{\lambda}y),\quad
 -\sqrt{\lambda}u_{x}\sin(\sqrt{\lambda}y),\quad\text{and}
\quad
 -\lambda u\cos(\sqrt{\lambda}y).
$$
 This fact allows us to estimate
the right-hand side of \eqref{06.5.18.5} and yields

\begin{equation}
                                                 \label{06.5.18.2}
I
\leq N\int_{Q_{2}}(|u_{xx}|^{p} +|u_{t}|^{p} 
+\lambda^{p/2}|u_{x}|^{p}+\lambda^{p}|u|^{p})\,dxdt.
\end{equation}
This is all we need since $\lambda|u|=|\bar{L}u |$ in $Q_{2}$ and the term
$\lambda^{p}|u|^{p}$ can be absorbed in $|u_{xx}|^{p} +|u_{t}|^{p} $.
 The lemma is proved.

Now comes the estimate of $v$ we were talking about
in Remark \ref{remark 06.10.21.1}. 

\begin{theorem}
                                        \label{theorem 06.5.15.1}
Let $\lambda\geq0$, $\nu\geq2$, 
and  $r\in(0,\infty)$ be some constants. Let  
 $u\in
C^{\infty}_{loc}(\bR^{d+1})$ be such that
$f:= \bar{L}u -\lambda u$ vanishes in $Q_{\nu r}$. Then
there is a constant $N=N(d,\kappa,K,p)$ such that
\begin{equation}
                                                  \label{06.5.15.1}
 (|u_{xx}(t,x) 
-(u_{xx})_{Q_{r}}|^{p} )_{Q_{r}}\leq N\nu^{-p}
 (|u_{xx}|^{p} +|u_{t}|^{p}
+\lambda^{p/2} |u_{x}|^{p} )_{Q_{\nu r}}.
\end{equation}

 \end{theorem}

Proof. Notice that $v(t,x):=u(tr^{2},xr)$ 
satisfy
$$
 (|u_{xx}(t,x) 
-(u_{xx})_{Q_{r}}|^{p} )_{Q_{r}}=r^{-2p}
(|v_{xx}(t,x) 
-(v_{xx})_{Q_{1}}|^{p} )_{Q_{1}},
$$
$$
 (|u_{xx}|^{p} +|u_{t}|^{p}
+\lambda^{p/2} |u_{x}|^{p} )_{Q_{\nu r}}
=r^{-2p} (|v_{xx}|^{p} +|v_{t}|^{p}
+\lambda^{p/2} r^{p}|v_{x}|^{p} )_{Q_{\nu}},
$$
and 
$$
\bar{L}(tr^{2})v(t,x) -r^{2}\lambda v(t,x)=
r^{2}f(tr^{2},xr)
$$ 
which vanishes in $Q_{\nu}$. It follows that if 
\eqref{06.5.15.1} holds for $r=1$, then it holds for any
$r>0$.

Therefore, in the rest of the proof we assume that $r=1$ and
observe that   the
left-hand side of \eqref{06.5.15.1} with $r=1$ is obviously
less than
a constant $N=N(d)$ times  
$$
\max_{Q_{1}}(|u_{xxx}|^{p}+|u_{txx}|^{p}).
$$

Therefore, we need only prove that
\begin{equation}
                                                  \label{06.5.15.2}
\max_{Q_{1}}(|u_{xxx}|^{p}+|u_{txx}|^{p})\leq
 N\nu^{-p}
 (|u_{xx}|^{p} +|u_{t}|^{p}
+\lambda^{p/2}|u_{x}|^{p} )_{Q_{\nu }}.
\end{equation}

Observe that the function $w(t,x)=u(t\nu^{2}/4,x\nu/2)$
satisfies
$$
\bar{L}(t\nu^{2}/4)w(t,x) -w(t,x)\nu^{2}\lambda /4=
0
$$
in $Q_{2}$ and
$$
 (|u_{xx}|^{p} +|u_{t}|^{p} +
\lambda^{p/2}|u_{x}|^{p})_{Q_{\nu }}
$$
$$
=(2/\nu)^{2p} (|w_{xx}|^{p} +|w_{t}|^{p} 
+(\nu^{2}\lambda/4)^{p/2}|w_{x}|^{p})_{Q_{2}},
$$
$$
 \max_{Q_{1}}|u_{xxx}|^{p}=(2/\nu)^{3p}
 \max_{Q_{2/\nu}}|w_{xxx}|^{p}\leq
(2/\nu)^{3p}
 \max_{Q_{1}}|w_{xxx}|^{p},
$$
$$
 \max_{Q_{1}} |u_{txx}|^{p}\leq(2/\nu)^{4p}
\max_{Q_{1}} |w_{txx}|^{p}.
$$

It follows that if \eqref{06.5.15.2} is true with $\nu=2$,
then
$$
\max_{Q_{1}}(|u_{xxx}|^{p}+|u_{txx}|^{p})\leq
N\nu^{-3p}\max_{Q_{1}}(|w_{xxx}|^{p}+|w_{txx}|^{p})
$$
$$
\leq N\nu^{-3p} (|w_{xx}|^{p} +|w_{t}|^{p} 
+(\nu^{2}\lambda/4)^{p/2}|w_{x}|^{p})_{Q_{2}}
$$
$$
= N\nu^{-p} (|u_{xx}|^{p} +|u_{t}|^{p} 
+\lambda^{p/2}|u_{x}|^{p})_{Q_{\nu }}.
$$
Finally, \eqref{06.5.15.2} with $\nu=2$ is indeed true 
by Lemma \ref{lemma 06.5.18.1} and the theorem is proved.

\begin{remark}
According to Theorem \ref{theorem 06.10.3.1},
applied to $u_{x}$ in place of $u$, the term  $|u_{t}|^{p}$
in \eqref{06.5.15.1} can be dropped.
\end{remark}

{\bf Proof of Theorem \ref{theorem 06.6.13.1}}. 
In Remark \ref{remark 06.10.21.1} we explained that
 we need $\lambda>0$ to
guarantee that certain equations have solutions.
 Therefore we take a $\lambda>0$, which in
the end will be sent to 0.  

Fix $r\in(0,\infty)$  and $\nu\geq4$.
We may certainly assume that
$a^{ij}$ are infinitely differentiable and have bounded
derivatives. Also changing $u$ for large $|t|+|x|$
does not affect \eqref{06.6.13.10}. Therefore, we may assume that
$u\in W^{1,2}_{p} $ and moreover $u\in
C^{\infty}_{0}(\bR^{d+1})$. In that case   define
$$
f=f_{\lambda}=\bar{L}u -\lambda u.
$$
Observe that $f\in C^{\infty}_{0}(\bR^{d+1})$.
Also take a $\zeta\in C^{\infty}_{0}(\bR^{d+1})$ such that
$\zeta=1$ on $Q_{\nu r/2}$ and $\zeta=0$ outside $Q_{\nu r}-Q_{\nu r}$
and set
$$
g=f\zeta,\quad h=f(1-\zeta).
$$
Finally define $v$ as the unique solution in $W^{1,2}_{p} $
of the equation
$$
\bar{L}v -\lambda v=h.
$$
Since $\lambda>0$, by classical theory we know that such a $v$ 
indeed exists
and is unique and infinitely differentiable. Since $h=0$ in
$Q_{\nu r/2}$ and $\nu/2\geq2$, by Theorem
\ref{theorem 06.5.15.1} we obtain
$$
(|v_{xx} 
-(v_{xx})_{Q_{r}}|^{p})_{Q_{r}}\leq N\nu^{-p}
(|v_{xx}|^{p} +|v_{t}|^{p}
+\lambda^{p/2} |v_{x}|^{p} )_{Q_{\nu r/2}}
$$
\begin{equation}
                                                  \label{06.6.13.3}
\leq N\nu^{-p}
 (|v_{xx}|^{p} +|v_{t}|^{p}
+\lambda^{p/2} |v_{x}|^{p} )_{Q_{\nu r}}.
\end{equation}

On the other hand the function 
$w:=u-v\in W^{1,2}_{p} $ satisfies
$$
\bar{L}w -\lambda w=g
$$
and by Theorem \ref{theorem 06.10.6.2}
$$
\int_{\bR^{d+1}_{0}}(|w_{t}|^{p}+|w_{xx}|^{p}+\lambda
^{p/2}|w_{x}|^{p})\,dxdt
$$
\begin{equation}
                                                  \label{06.6.13.4}
\leq N\int_{\bR^{d+1}_{0}}|g|^{p}\,dxdt
\leq N\int_{Q_{\nu r}}|f|^{p}\,dxdt,
\end{equation}
$$
\int_{Q_{r}}|w_{xx}|^{p}\,dxdt
\leq N\int_{Q_{\nu r}}|f|^{p}\,dxdt,
$$
\begin{equation}
                                                  \label{06.6.14.3}
(|w_{xx}|^{p})_{Q_{r}}
\leq N\nu^{d+2}(|f|^{p})_{Q_{\nu r}}.
\end{equation}
By combining this with \eqref{06.6.13.3} and observing that $u=v+w$ and
$$
I:=(|u_{xx} 
-(u_{xx})_{Q_{r}}|^{p})_{Q_{r}}
\leq2^{p}(|w_{xx} 
-(w_{xx})_{Q_{r}}|^{p})_{Q_{r}}
$$
$$
+2^{p}(|v_{xx}
-(v_{xx})_{Q_{r}}|^{p})_{Q_{r}}
\leq N(|w_{xx}|^{p} )_{Q_{r}}+
2^{p}(|v_{xx} 
-(v_{xx})_{Q_{r}}|^{p} )_{Q_{r}},
$$
we get
$$
I\leq N\nu^{d+2}(|f|^{p})_{Q_{\nu r}}
+N\nu^{-p}
 (|v_{xx}|^{p} +|v_{t}|^{p}
+\lambda^{p/2}|v_{x}|^{p} )_{Q_{\nu r}} 
$$
$$
\leq
N\nu^{d+2}(|f|^{p})_{Q_{\nu r}}
+N\nu^{-p}
 (|u_{xx}|^{p} +|u_{t}|^{p}
+\lambda^{p/2} |u_{x}|^{p} )_{Q_{\nu r}}
$$
$$
+N\nu^{-p}
 (|w_{xx}|^{p} +|w_{t}|^{p}
+\lambda^{p/2} |w_{x}|^{p} )_{Q_{\nu r}}.
$$
Here by \eqref{06.6.13.4}
$$
 (|w_{xx}|^{p} +|w_{t}|^{p}
+\lambda^{p/2} |w_{x}|^{p} )_{Q_{\nu r}}\leq
 N(|f|^{p})_{Q_{\nu r}}
$$
and since $\nu\geq1$ we conclude
$$
I\leq
N\nu^{d+2}(|f_{\lambda}|^{p})_{Q_{\nu r}}
+N\nu^{-p}
 (|u_{xx}|^{p} +|u_{t}|^{p}
+\lambda^{p/2} |u_{x}|^{p} )_{Q_{\nu r}}.
$$
To get \eqref{06.6.13.10} it only remains to
use that $u_{t}=f_{\lambda}+\lambda
u-a^{ij}u_{x^{i}x^{j}}$
and let $\lambda
\downarrow0$. The theorem is proved.
\begin{remark}
                                          \label{remark 06.10.17.3}
Recall that for $\phi\in L_{1,loc}$ the sharp
function $\phi^{\#}$ and the maximal function $M\phi$ are defined by
$$
\phi^{\#}(t,x)=\sup_{Q\in\bQ:(t,x)\in Q}
(|\phi-\phi_{Q}|)_{Q},\quad Mf(t,x)=
\sup_{Q\in\bQ:(t,x)\in Q}
\phi_{Q}.
$$
In this notation Theorem \ref{theorem 06.6.13.1} 
and H\"older's inequality imply that
on $\bR^{d+1}$ we have
$$
(u_{xx})^{\#}\leq N\nu^{(d+2)/p}M^{1/p}(|\bar{L}u|^{p})
+N\nu^{-1}M^{1/p}(|u_{xx}|^{p}).
$$
Then by using the Fefferman-Stein theorem we obtain for any $q>p$
$$
\|u_{xx}\|_{L_{q}}\leq N\|(u_{xx})^{\#}\|_{L_{q}}
\leq N\nu^{(d+2)/p}\|\bar{L}u\|_{L_{q}}
+N\nu^{-1}\|u_{xx}\|_{L_{q}},
$$
where the second inequality holds
since $\|M^{1/p}\phi\|_{L_{q}}\leq N\| \phi^{1/p}\|_{L_{q}}$
by the Hardy-Littlewood theorem. For $\nu$ large enough
we absorb the last term into the left-hand side and get
\begin{equation}
                                                      \label{06.10.17.3}
\|u_{xx}\|_{L_{q}}\leq N \|\bar{L}u\|_{L_{q}}.
\end{equation}

This and what is said after Theorem \ref{theorem 06.6.13.1}
allow us to give one more proof of Theorem \ref{theorem 06.10.6.2}.
 
\end{remark}  

To summarize,
after having proved  Theorem \ref{theorem 06.6.13.1} 
one can follow the
same way as in \cite{Kr06} and get the solvability of equations with
$VMO_{x}$ leading coefficients.

In particular, we have the following result.  

\begin{theorem}
                                                   \label{theorem 7.14.1}
There are   constants
$\lambda_{0}$ and
$N$, depending only on $p$, $K$, $\kappa$, $d$, and 
$\omega$, such that for any $\lambda\geq\lambda_{0}$
and $u\in W^{1,2}_{p}$ we have  
\begin{equation}
                                                     \label{7.14.3}
\lambda\|u\|_{L_{p}}+\sqrt{\lambda}\|u_{x }\|_{L_{p}}
+\|u_{xx}\|_{L_{p}}+\|u_{t}\|_{L_{p}}\leq N\|(L-\lambda)u\|_{L_{p}}.
\end{equation}
 
Furthermore, for any $\lambda\geq\lambda_{0}$ and
$f\in L_{p}$ there exists a unique 
$u\in W^{1,2}_{p}$ such that $(L-\lambda)u=f$.

\end{theorem}  

\begin{corollary}
                                            \label{corollary 06.10.5.1}
There is a constant 
$N_{0}$, depending only on $p$, $K$, $\kappa$, $d$, and 
$\omega$, such that for any   $u\in W^{1,2}_{p}(\bR^{d+1}_{0})$ we have  
\begin{equation}
                                                     \label{06.10.5.2}
 \|u_{xx}\|_{L_{p}(\bR^{d+1}_{0})}
+\|u_{t}\|_{L_{p}(\bR^{d+1}_{0})}\leq
N_{0}(\| L u\|_{L_{p}(\bR^{d+1}_{0})}+
\|  u\|_{L_{p}(\bR^{d+1}_{0})}).
\end{equation}
\end{corollary}

To prove this we first claim that \eqref{7.14.3}
with $L_{p}(\bR^{d+1}_{0})$ in place of $L_{p}$
holds for any $u\in W^{1,2}_{p}(\bR^{d+1}_{0})$.

Indeed, for such a $u$ set $f(t,x)=I_{t>0}(L-\lambda)u(t,x)$,
let $v\in W^{1,2}_{p}$ be any function
on $\bR^{d+1}$ coinciding with $u$ for $t>0$, and set
$g=(L-\lambda)v$. Then find $w\in W^{1,2}_{p}$ such that
$(L-\lambda)w=f$ and observe that $(L-\lambda)(v-w)=g-f$
vanishes for $t>0$. One can solve the equation
$(L-\lambda)\phi=g-f$ by the method of continuity
starting from $L=\Delta+D_{t}$,
for which the solutions vanish for $t>0$
if the right-hand side does that, and then one sees
that $v=w$ for $t>0$. This means that $u=w$
for $t>0$. Since estimate \eqref{7.14.3} holds
with $w$ in place of $u$ and $f$ in place of $(L-\lambda)u$,
we get our claim.

After that it suffices to take $\lambda=\lambda_{0}$ and observe that
$$
\| L u-\lambda_{0} u\|_{L_{p}(\bR^{d+1}_{0})}
\leq \| L u\|_{L_{p}(\bR^{d+1}_{0})}
+\lambda_{0}\|   u\|_{L_{p}(\bR^{d+1}_{0})}.
$$

\mysection{Proof of Lemma \protect\ref{lemma 06.7.26.2}} 
                                      \label{section 06.10.19.1}

The program of proof is to use Theorem \ref{theorem 06.6.13.1}
but replace $\bar{L}u$ in \eqref{06.6.13.10}
with $Lu$. The error term we estimate by using
H\"older's inequality and on the account of
right choice of $\bar{L}$ come to \eqref{06.7.26.8}
with
$$
(|u_{xx}|^{2p})_{Q_{\rho}}^{1/2}
$$
 in place of $\cB_{\rho}$. Then the main issue is how to
reduce power $2p$ back to $p$.
It turns out that this is possible if $u$ is ``harmonic"
in $Q_{2\rho}$ (see Corollary \ref{corollary 06.7.7.1}). After that
 we use the same kind of decomposition
of $u$ as in Remark \ref{remark 06.10.21.1}.
As in Lemma \ref{lemma 06.7.26.2} we assume that 
$p\in(1,\infty)$, $b=0$, and $c=0$.

We   need two versions of Lemma \ref{lemma 06.7.7.2}
when the powers of summability on the right are less
than on the left. Similar estimate  is known even with $\nu=1$
for the elliptic case as Poincar\'e's inequality.
\begin{lemma}
                                          \label{lemma 06.7.26.1}
Let $q\geq1$, $\nu\in(1,\infty)$,
\begin{equation}
                                                 \label{06.7.26.1}
 \frac{1}{q}<\frac{2}{d+2}+\frac{1}{p}.
\end{equation}
Then there is a constant $N=N(d,p,q,\nu)$ such that
for any $u\in W^{1,2}_{q,loc}$ and
$r\in(0,\infty)$ we have
\begin{equation}
                                                 \label{06.7.26.5}
 ( |u(t,x)- u _{Q_{\nu r}}-
x^{i}(u_{x^{i}})_{Q_{\nu r}}|^{p} )^{1/p}_{Q_{r}}
\leq Nr^{2}(  |u_{xx}|^{q} +|u_{t}|^{q}  
 )^{1/q}_{Q_{\nu r}}.
\end{equation}
\end{lemma}

Proof.  First, observe that an argument based on
 self-similarity reduces the case of general $r$ to the case that $r=1$,
the one we confine ourselves to.
Then by obvious reasons we may assume that $u\in
C^{\infty}_{0}(\bR^{d+1})$. Finally, if $q\geq p$, the
result follows from Lemma \ref{lemma 06.7.7.2} and H\"older's
inequality. Therefore, we assume that $q\leq p$.

Take an infinitely differentiable function $\zeta$ on $\bR^{d+1}$
such that $\zeta=1$ on $Q_{1}$ and $\zeta=0$ on $\bR^{d+1}_{0}
\setminus Q_{\nu}$,
and set  
$$
f= \Delta u+u_{t} ,\quad v=\zeta 
(u-u_{Q_{\nu}}-x^{i}(u_{x^{i}})_{Q_{\nu}}),
$$
so that
$$
\Delta  v+ v_{t} 
=\zeta f+(u-u_{Q_{\nu}}-x^{i}(u_{x^{i}})_{Q_{\nu}})
(\Delta\zeta+\zeta_{t})+2\zeta_{x^{i}}(u_{x^{i}}-(u_{x^{i}})_{Q_{\nu}})
=:-g.
$$

Since $v\in C^{\infty}_{0}(\bR^{d+1})$, we have
$$
 v(t,x)=\int_{0}^{\infty}\int_{\bR^{d}}
g(t+s,x+y) p(s,y)\,dyds,\quad p(s,y)=\frac{1}{(4\pi s)^{d/2}}
e^{-|y|^{2}/(4s)}.
$$ 
Here, if $0\leq t\leq 1$, there is no need to integrate
with respect to $s$ beyond $[0,\nu^{2}]$, since $g(r,z)=0$ for $r\geq
\nu^{2}$. Therefore, upon denoting
$$
\bar{v}(t,x)=|v (t,x)|I_{t\in[0,1]},\quad
\bar{g}(s,y)=|g(s,y)|I_{s\in[0,\nu^{2}]},
$$
$$
\bar{p}(s,y)=p(s,y)I_{s\in[0,\nu^{2}]},
$$
we find
$$
\bar{v}(t,x)\leq  \int_{\bR^{d+1}}\bar{g}(t+s,x+y)\bar{p}(s,y)\,dxds
$$
$$
= \int_{\bR^{d+1}}\bar{g}(t-s,x-y)\bar{p}(-s,y)\,dxds.  
$$

Now we  apply Young's inequality 
\begin{equation}
                                              \label{06.10.7.1}
\|\bar{g}*\bar{p}\|_{L_{p} }
\leq\|\bar{g}\|_{L_{q}  }\|\bar{p}\|_{L_{r} },
\end{equation}
where
$$
r=\frac{pq}{q-p+p q},
$$
  $r\geq1$ since $q\leq p$, and $p^{-1}+1=q^{-1}+r^{-1}$.
Also 
\begin{equation}
                                                 \label{06.7.26.2}
rd< d+2 
\end{equation}
due  to \eqref{06.7.26.1}. Then we find
\begin{equation}
                                                 \label{06.7.26.4}
\|v\|_{L_{p}(Q_{1})}\leq\|\bar{v}\|_{L_{p} }
\leq  \|g\|_{L_{q}(Q_{\nu})}\|p
\|_{L_{r}([0,\nu^{2}]\times\bR^{d})}.
\end{equation}

Here by the definition of $g$ and Lemma \ref{lemma 06.7.7.2}
(just in case, recall that $N$ in \eqref{06.7.26.5}
is allowed to depend on $\nu$)
$$
\|g\|_{L_{q}(Q_{\nu})}
\leq N(\|u_{xx}\|_{L_{q}(Q_{\nu})}+\|u_{t}\|_{L_{q}(Q_{\nu})}
 ).
$$
Furthermore, changing variables shows that the integral
$$
\int_{\bR^{d}}t^{-d/2}
e^{-r|x|^{2}/(4t)}\,dx
$$
is finite and independent of $t>0$. Therefore,
$$
\|p\|_{L_{r}([0,4]\times\bR^{d})}^{r}
=N\int_{0}^{4}t^{-rd/2+d/2}\int_{\bR^{d}}t^{-d/2}
e^{-r|x|^{2}/(4t)}\,dxdt
$$
$$
=N\int_{0}^{4}t^{-rd/2+d/2}\,dt<\infty,
$$
where the inequality holds since
owing to \eqref{06.7.26.2} we have $-rd/2+d/2>-1$.

Thus, \eqref{06.7.26.4} implies that
$$
\|v\|_{L_{p}(Q_{1})}\leq N(\|u_{xx}\|_{L_{q}(Q_{\nu})}+
\|u_{t}\|_{L_{q}(Q_{\nu})} )
$$
and it only remains to observe that the left-hand side here
coincides with the left-hand side of \eqref{06.7.26.5}.
 The lemma is proved.

Similar estimate holds for $u_{x}-(u_{x})_{Q_{\nu r}}$.
\begin{lemma}
                                            \label{lemma 06.7.17.1} 
Let $q\geq1$, $\nu\in(1,\infty)$,
\begin{equation}
                                                 \label{06.7.18.3}
\frac{1}{q}<\frac{1}{d+2}+\frac{1}{p} .
\end{equation}
Then there is a constant $N=N(d,p,q,\nu)$ such that
for any $u\in W^{1,2}_{q,loc}$ and
$r\in(0,\infty)$ we have
\begin{equation}
                                                 \label{06.10.7.2}
 ( |u_{x}(t,x) 
-(u_{x})_{Q_{r}}|^{p}  )^{1/p}_{Q_{r}}
\leq Nr (  |u_{xx}|^{q} +|u_{t}|^{q}  
 )^{1/q}_{Q_{\nu r}}.
\end{equation}
\end{lemma}

Proof. As in the proof of Lemma \ref{lemma 06.7.26.1}
we may assume that $r=1$, $q\leq p$,
 and $u\in C^{\infty}_{0}(\bR^{d+1})$.
Then, again
take an infinitely differentiable function $\zeta$ on $\bR^{d+1}$
such that $\zeta=1$ on $Q_{1}$ and $\zeta=0$ on $\bR^{d+1}_{0}
\setminus Q_{\nu}$, and use the notation from the proof
of Lemma \ref{lemma 06.7.26.1} to obtain
$$
v_{x}(t,x)=\int_{0}^{\infty}\int_{\bR^{d}}
g(t+s,x+y)p_{y}(s,y)\,dyds
$$

Next, we use an elementary inequality
$$
x^{\alpha}e^{-\beta x}\leq Ne^{-\beta x/2},\quad\forall x\geq0,
$$
where $\alpha,\beta>0$ and $N=N(\alpha,\beta)$. Then by observing that
$$
p_{y^{i}}(s,y)=-\frac{y^{i}}{2s}\frac{1}{(4\pi s)^{d/2}}e^{-|y|^{2}/(4s)}
$$
we find
$$
|p_{y}(s,y)|\leq
\frac{1}{\sqrt{s}}
\frac{|y|}{\sqrt{4s}}\frac{1}{(4\pi s)^{d/2}}e^{-|y|^{2}/(4s)}
\leq Ns^{-1/2}p(s/2,y),
$$
which implies that
$$
| v_{x}(t,x)|\leq N
\int_{0}^{\infty}\int_{\bR^{d}}|g(t+s,x+y)|s^{-1/2}
 p(s/2,y)\,dyds.
$$
As before,
 if $0\leq t\leq 1$, there is no need to integrate
with respect to $s$ beyond $[0,\nu^{2}]$.
Therefore, upon denoting
$$
w(t,x)=|v_{x}(t,x)|I_{t\in[0,1]},\quad
\bar{g}(s,y)=|g(s,y)|I_{s\in[0,\nu^{2}]},
$$
$$
 h(s,y)=s^{-1/2} p(s/2,y)
I_{s\in[0,\nu^{2}]},
$$
we find
$$
w(t,x)\leq N\int_{\bR^{d+1}}\bar{g}(t+s,x+y)h(s,y)\,dxds
$$
$$
=N\int_{\bR^{d+1}}\bar{g}(t-s,x-y)h(-s,y)\,dxds.
$$

Now we  apply \eqref{06.10.7.1}
with the same $r$, which also satisfies
\begin{equation}
                                                 \label{06.7.18.2}
r(d+1)< d+2 
\end{equation}
due  to \eqref{06.7.18.3}. Then we find
\begin{equation}
                                                 \label{06.7.18.4}
\|v_{x}\|_{L_{p}(Q_{1})}\leq\|w\|_{L_{p} }
\leq N\|g\|_{L_{q}(Q_{\nu})}\|h\|_{L_{r}([0,\nu^{2}]\times\bR^{d})}.
\end{equation}

Here by the definition of $g$ and Lemma \ref{lemma 06.7.7.2}
$$
\|g\|_{L_{q}(Q_{\nu})}\leq 
N(\|u_{xx}\|_{L_{q}(Q_{\nu})}+\|u_{t}\|_{L_{q}(Q_{\nu})}).
$$
Furthermore, 
$$
\|h\|_{L_{r}([0,4]\times\bR^{d})}^{r}
=N\int_{0}^{4}t^{-r(d+1)/2+d/2}\int_{\bR^{d}}t^{-d/2}
e^{-r|x|^{2}/(2t)}\,dxdt
$$
$$
=N\int_{0}^{4}t^{-r(d+1)/2+d/2}\,dt<\infty,
$$
where the inequality holds since
owing to \eqref{06.7.18.2} we have $-r(d+1)/2+d/2>-1$.

Now it only remains to observe that the
left-hand sides of \eqref{06.7.18.4}
and \eqref{06.10.7.2} coincide.
 The lemma is proved.
\begin{lemma}
                                 \label{lemma 06.10.5.1}
 Let $r\in(0,1]$, $\nu\in(1,\infty)$, and
$u\in W^{1,2}_{p,loc}$.
Set
$f:=Lu $.
 Then  
\begin{equation}
                                               \label{06.7.7.7}
 \|u_{xx}\|_{L_{p}(Q_{r})} 
\leq N\big(
\|f\|_{L_{p}(Q_{\nu r})} 
+r^{-1} \|u_{x}\|_{L_{p}(Q_{\nu r})} 
+r^{-2} \|u\|_{L_{p}(Q_{\nu r})}\big),
\end{equation}
where $N$ depends only on $\nu, d,K,p, \kappa$,
and the function $\omega$. 
\end{lemma}

Proof. Obviously we may concentrate on 
$u\in W^{1,2}_{p }(\bR^{d+1}_{0})$.
By Corollary \ref{corollary 06.10.5.1}
the assumption of Lemma \ref{lemma 9.15.1}
is satisfied. Therefore, \eqref{06.7.7.7}  
holds with $r=1$.

For $r\in(0,1]$ and $u\in W^{1,2}_{p }(\bR^{d+1}_{0})$
introduce $v(t,x)=u(r^{ 2}t,rx  )$ and observe that
$$
v_{t}(t,x)+\bar{a}^{ij}(t,x)v_{x^{i}x^{j}}(t,x)
+rb^{i}(r^{ 2}t,rx )v_{x^{i}}(t,x)
$$
\begin{equation}
                                                  \label{06.10.5.3}
+r^{2}c(r^{ 2}t,rx )
v(t,x)=r^{ 2}f(r^{ 2}t,rx )=:g(t,x),
\end{equation}
where $\bar{a}(t,x)=a (r^{ 2}t,rx )$.

Furthermore, for any $\rho>0$ and $t,x$
$$
\rho^{-2}|B_{\rho}|
\int_{t}^{t+\rho^{2}}\int_{y,z\in B_{\rho}(x)}
|\bar{a}(s ,y )-\bar{a}(s ,z )|\,dyds 
$$
$$
=(r\rho)^{-2}|B_{r\rho}|
\int_{r^{2}t}^{r^{2}t+(r\rho)^{2}}\int_{y,z\in B_{r\rho}(rx)}
|a (s ,y )-a (s ,z )|\,dyds .
$$
Therefore, $\bar{a}_{\rho}^{\#(x)}\leq\omega(r\rho)\leq
\omega(\rho)$. Also $|rb|\leq K$ and $r^{2}|c|\leq K$.
It follows that the above result is applicable
to \eqref{06.10.5.3} and
$$
 \|v_{xx}\|_{L_{p}(Q_{1})} 
\leq N\big(
\|g\|_{L_{p}(Q_{2 })}  
+  \|v_{x}\|_{L_{p}(Q_{2 })}
+  \|v\|_{L_{p}(Q_{2 })}\big).
$$
Expressing all terms here by means of $u$ and $f$ leads to 
\eqref{06.7.7.7}. The lemma is proved.

The following is a  crucial 
point in proving Lemma  \ref{lemma 06.7.26.2}.
 
\begin{corollary}
                                           \label{corollary 06.7.7.1}
If
  $r\in(0,1]$, $q\geq1$, and
$u\in W^{1,2}_{p,loc} $ 
are such that that in $Q_{2 r}$ we have
$Lu  =0$, $b=0$, and $c=0$,
then
\begin{equation}
                                                  \label{06.7.30.2}
(|u_{xx}|^{p} )^{1/p}_{Q_{r}}
\leq N_{1} (|u_{xx}|)_{Q_{2r}} 
\leq
N_{1} (|u_{xx}|^{q}  )^{1/q}_{Q_{2r}} ,
\end{equation}
where   
 $N_{1}$ depends only on $ d ,p,
\kappa,K
$, and the function $\omega$.
\end{corollary}

Proof. The second inequality  in \eqref{06.7.30.2}
follows from H\"older's inequality. It turns out that,
to prove the first one, it
suffices to prove that if \eqref{06.7.18.3} holds, $q\leq p$,
$\nu\in(1,\infty)$,
and $Lu=0$ in $Q_{\nu r}$,
then
\begin{equation}
                                                  \label{06.7.7.5}
(|u_{xx}|^{p} )^{1/p}_{Q_{r}}
\leq N(|u_{xx}|^{q} )^{1/q}_{Q_{\nu r}},
\end{equation}
where   $N=N(\nu,d, p,q,\omega,\kappa,K )$.
Indeed, one can find a decreasing sequence $q_{i}\in[1,p]$,
$i=0,1,...,m$, where $m$ depends only on $p$
and $d$, such that $q_{0}=p$, $q_{m}=1$, and 
$q_{i+1}^{-1}<(d+2)^{-1}+q_{i}^{-1}$.
Then if  \eqref{06.7.7.5} is true
under the additional assumptions, then 
the $L_{q_{i}}$ average norm of $u_{xx}$ is estimated
by the $L_{q_{i+1}}$ average norm of $u_{xx}$  
in an expanded domain of averaging.
We can then iterate \eqref{06.7.7.5} going 
along the sequence $q_{i}$ and
we can choose $\nu=\nu(p)$ so close to 1, that during these
 finitely many
steps the expanding domains would always be in $Q_{2r}$
and 
\eqref{06.7.30.2} would follow.

Therefore, we concentrate on proving \eqref{06.7.7.5} assuming that
 \eqref{06.7.18.3} holds, $q\leq p$,
$\nu\in(1,\infty)$,
and $Lu=0$ in $Q_{\nu r}$
 Since \eqref{06.7.7.5} only involves the values of $u$
in $Q_{\nu r}$, we may assume that 
$u\in W^{1,2}_{p } $. In that case introduce
$$
v = u-u_{Q_{\nu r}}-x^{i}(u_{x^{i}})_{Q_{\nu r}}.
$$
Since by assumption $Lv  =0$
  in $Q_{\nu r}$ and $r\leq 1$,
by Lemma \ref{lemma 06.10.5.1}  
  we have
$$
\int_{Q_{r}}|u_{xx}|^{p}\,dxdt=\int_{Q_{r}}|v_{xx}|^{p}\,dxdt
\leq 
Nr^{- p}
\int_{Q_{\sqrt{\nu} r}}|u_{x} -
 (u_{x })_{Q_{\nu r}}|^{p}\,dxdt
$$
\begin{equation}
                                                 \label{06.7.26.7}
+
Nr^{-2p}
\int_{Q_{\sqrt{\nu} r}}|u-u_{Q_{\nu r}}-
x^{i}(u_{x^{i}})_{Q_{\nu r}}|^{p}\,dxdt.
\end{equation}
By Lemmas \ref{lemma 06.7.26.1}  and \ref{lemma 06.7.17.1} 
the right-hand side in \eqref{06.7.26.7} is less than
the $p$-th power of the right-hand side in \eqref{06.7.7.5}.
 The corollary is proved.

{\bf Proof of Lemma  \ref{lemma 06.7.26.2}}.
  According to Theorem 2.1 of \cite{Kr06}, 
 there is a function $v$ such that
it belongs to $W^{1,2}_{p}(\bR^{d+1}_{S})$ for any $S>-\infty$,
satisfies
\begin{equation}
                                                \label{06.7.28.1}
Lv =fI_{Q_{\nu r}} 
\end{equation}
in $\bR^{d+1} $, and is such that $v(t,x)=0$ for $t>4$
(observe that $\nu r\leq1$). 
Furthermore, as usual, since $fI_{Q_{\nu r}}\in L_{q} $
for any  $q\in(1,\infty)$,
we have that
$v\in  W ^{1,2}_{q}(\bR^{d+1}_{S})$ for all   $q\in(1,\infty)$
and $S $.

After that we  set
$$
w=u-v
$$ 
and note for the future that
$w\in W ^{1,2}_{q,loc} $ for all  $q\in(1,\infty)$.

Again by Theorem 2.1 of \cite{Kr06} we have
$$
\int_{(0,4)\times\bR^{d}}
|v_{xx}|^{p}\,dxdt\leq N\int_{Q_{\nu r}}|f|^{p}\,dxdt
$$
implying that
\begin{equation}
                                                \label{06.7.28.2}
(|v_{xx}|^{p})_{Q_{\nu r}}\leq N\cA_{\nu r},\quad
(|v_{xx}|^{p})_{Q_{r}}\leq N\nu^{d+2}\cA_{\nu r}.
\end{equation}

Next, observe that 
$$
w\in W^{1,2}_{2p,loc } 
\subset W^{1,2}_{p,loc} 
$$
and $Lw =0$ in $Q_{\nu r}$ and $\nu/4\geq4$. 

Now we apply Theorem \ref{theorem 06.6.13.1}
with $\nu/4$ in place of $\nu$,
$\bar{L}w=w_{t}+\bar{a}^{ij}w_{x^{i}x^{j}}$ and $\bar{a}\in\bA$. As an
intermediate step we also use H\"older's inequality and the fact that
$\bar{L}w=(\bar{a}-a )^{ij} w_{x^{i}x^{j}}$ in $Q_{\nu r}$ to find that  
$$
\dashint_{Q_{\nu r/4}}|\bar{L}w|^{p}\,dxdt
\leq N
\big(\dashint_{Q_{\nu r/4}}|w_{xx}|^{2p}\,dxdt\big)^{1/2}
\big(\dashint_{Q_{\nu r/4}}|a-\bar{a} |^{2p}\,dxdt\big)^{1/2},
$$  
where for for an appropriate
$\bar{a}$
$$
\big(\dashint_{Q_{\nu r/4}}|a-\bar{a} |^{2p}\,dxdt\big)^{1/2}
\leq
N \big(\dashint_{Q_{\nu r/4}}|a-\bar{a} | \,dxdt\big)^{1/2}
\leq N\hat{a}^{1/2}.
$$

Then we obtain
\begin{equation}
                                               \label{06.7.28.3}
 (|w_{xx}-(w_{xx})_{Q_{r} }|^{p})_{Q_{r} }
\leq  N\nu^{-p}(|w_{xx}|^{p})_{Q_{\nu r/4}}+
N \nu ^{d+2} 
  \hat{a}^{1/2}[(|w_{xx}|^{2 p})_{Q_{\nu r/4}}]^{1/2}. 
\end{equation}

Owing to \eqref{06.7.28.2} and the definition of $w$,
$$
(|w_{xx}|^{p})_{Q_{\nu r/4}}\leq 
N(|w_{xx}|^{p})_{Q_{\nu r }}\leq
N(|w_{xx}+v_{xx}|^{p})_{Q_{\nu r }}
$$
\begin{equation}
                                               \label{06.7.28.5}
+N(|v_{xx}|^{p})_{Q_{\nu r }}
\leq N\cB_{\nu r}+N\cA_{\nu r}.
\end{equation}
Now we apply Corollary \ref{corollary 06.7.7.1}
with $2p$ in place of $p$  noting that
  the fact that
$Lw =0$ in $Q_{\nu r}$ allows us to do that.
Then we see that
$$
[(|w_{xx}|^{2 p})_{Q_{\nu r/4}}]^{1/2}
\leq N(|w_{xx}|^{p})_{Q_{\nu r}}.
$$
We estimate the last term using \eqref{06.7.28.5}
and then infer from \eqref{06.7.28.3} that
$$
(|w_{xx}-(w_{xx})_{Q_{r} }|^{p})_{Q_{r} }
\leq  N(\nu^{-p}+\nu ^{d+2} 
  \hat{a}^{1/2})(\cB_{\nu r}+ \cA_{\nu r}).
$$
To finish proving \eqref{06.7.26.8} it only remains
to combine this with \eqref{06.7.28.2} and observe that
$$
(|u_{xx}-(u_{xx})_{Q_{r} }|^{p})_{Q_{r} }\leq
N(|v_{xx}-(v_{xx})_{Q_{r} }|^{p})_{Q_{r} }
+N(|w_{xx}-(w_{xx})_{Q_{r} }|^{p})_{Q_{r} },
$$
$$
(|v_{xx}-(v_{xx})_{Q_{r} }|^{p})_{Q_{r} }\leq
N(|v_{xx} |^{p})_{Q_{r} }.
$$
The lemma is proved.

\mysection{New approach to the $L_{p}$-theory for
divergence type equations  
with VMO coefficients}  
                                       \label{section 7.22.20}

Take an $a\in\bA$ and set
$$
\bar{L}u(t,x)=a^{ij}(t)u_{x^{i}x^{j}}(t,x)+u_{t}(t,x).
$$
In this section we show how to use results on solvability
of equations in the whole space and prove
the following statement which is a weak version
of Lemma \ref{lemma 06.10.9.5} and
for $p=2$ is Lemma 5.2 of \cite{Kr06} proved there by using the
solvability of equations in cylinders.
Throughout the section
  $p\in(1,\infty)$ and $\lambda\geq0$
unless explicitly specified otherwise.

\begin{theorem}
                                        \label{theorem 06.10.6.1}
Let  
$u\in
\cH^{1}_{p,loc} $, $f=(f^{1},...,f^{d})$,
 $f^{i}\in L_{p,loc}$, $\nu\geq4$, $r>0$.
Assume that $\bar{L}u=\Div f$ in $Q_{\nu r}$. Then  
there exists a constant
$N=N( d,\kappa,K ,p)$ such that  
\begin{equation}
                                                 \label{06.10.6.1}
\big(|u_{x }-( u_{x })_{Q_{r}}| ^{p}\big)_{Q_{r}}
\leq N\nu^{-p}
\big(|u_{ x}|^{p} \big)_{Q_{\nu r}} 
+N \nu^{ d+2 }
\big(  |f|^{p}   \big)_{Q_{\nu r}} .
\end{equation}
 \end{theorem}

Our strategy is very similar to what is done
in Section \ref{section 7.22.2}.
We need few auxiliary results.
The first  one is   used also later in the proof of Corollary
\ref{corollary 06.10.8.3}.

\begin{lemma}
                                                 \label{lemma 7.8.1}
Let $p\in[1,\infty)$,
 $R\in(0,\infty)$,  $u\in
\cH^{1}_{p,loc}$,
$$
f=(f^{1},...,f^{d}),\quad  f^{i},g\in L_{p,loc} ,
$$
and  $\bar{L}u=\Div f+g$ in $Q_{R}$. Then
for a constant $N=N(d,K,p)$ we have
\begin{equation}
                                                     \label{7.8.1}
\int_{Q_{R}}|u(t,x)- u_{Q_{R}}|^{p} \,dxdt
\leq NR^{p} \int_{Q_{R}}(|u_{x }|^{p} +|f|^{p}+R^{p}|g|^{p} )\,dxdt.
\end{equation}

\end{lemma}

Proof. Denote by $\phi^{(\varepsilon)}$ the convolution
of $\varepsilon^{-d-2}\zeta(\varepsilon
^{-2}t,\varepsilon^{-1}x)$ with  $\phi=\phi(t,x)$, where
 $\zeta\in C^{\infty}_{0}(\bR^{d+1})$ and $\zeta$
integrates to one. Let $\bar{L}^{(\varepsilon)}$
be the operator constructed from $a^{(\varepsilon)}$.
Observe that
the equation
\begin{equation}
                                                 \label{06.10.10.05}
\bar{L}^{(\varepsilon)}u^{(\varepsilon)}=\Div 
f^{\varepsilon}+
g^{\varepsilon},
\end{equation}
where
$$
f^{\varepsilon j}=f^{(\varepsilon)j}+
  a^{(\varepsilon)ij}
u^{(\varepsilon)}_{x^{i}}-(a^{ij}u_{x^{i}})^{(\varepsilon)},
$$
holds in a somewhat smaller domain than $Q_{R}$.
If the assertion of the lemma were 
applicable to \eqref{06.10.10.05} and somewhat smaller domains, then,
since $u^{(\varepsilon)}$, $u^{(\varepsilon)}_{x}$, $f^{\varepsilon}$,
and
$g^{\varepsilon}$ converge in $L_{p} $
as $\varepsilon\to0$
to $u$, $u_{x}$, $f$, and $g$, respectively, we would get 
\eqref{7.8.1}.
This argument convinces us that without 
losing generality we may assume that
$a,u,f$, and $g$ are infinitely differentiable.
In that case our assertion is known as is Lemma 3.1 of
\cite{Kr06}. The lemma is proved.

\begin{lemma}
                                          \label{lemma 06.10.3.1}
Let   $m\in\{0,1,2,...\}$  and
 $u\in C^{\infty}_{0}(\bR^{d+1})$. Assume that 
$\bar{L}u -\lambda u$ vanishes in $Q_{2}$. Then
\begin{equation}
                                                 \label{06.10.3.1}
\max_{Q_{1}}\big(|D^{m}u_{x }|^{p}
+|D^{m}u_{t}|^{p}\big)
\leq N\int_{Q_{2}}(|u_{x }|^{p} 
+\lambda^{p/2}|u |^{p})\,dxdt,
\end{equation}
where $N=N(d,m,\kappa,K,p)$.
\end{lemma}

Proof.
If $\lambda=0$, we obtain the estimate of $|D^{m}u_{x}|$
by applying \eqref{06.5.27.3} with $|\alpha|\geq1$
to $u-u_{Q_{r}}$ in place of $u$ and 
using Lemma \ref{lemma 7.8.1}. The estimate for
$D^{\alpha}u_{t}$ then follows from the equation $\bar{L}u=0$
in $Q_{2}$.

For general $\lambda$ we just inspect the proof 
of Lemma \ref{lemma 06.5.18.1} and observe that
it works in the present case as well.
The lemma is proved.

Here is a counterpart of Theorem \ref{theorem 06.5.15.1}
which is proved in the same way.  
\begin{theorem}
                                        \label{theorem 06.10.3.1}
Let $\lambda\geq0$, $\nu\geq2$,
and  $r\in(0,\infty)$ be some constants. Let  
 $u\in
C^{\infty}_{loc}(\bR^{d+1})$ be such that
$ \bar{L}u -\lambda u$ vanishes in $Q_{\nu r}$. Then
there is a constant $N=N(d,\kappa,K,p)$ such that
\begin{equation}
                                                  \label{06.10.3.2}
(|u_{x } 
-(u_{x })_{Q_{r}}|^{p})_{Q_{r}}\leq N\nu^{-p}
(|u_{x }|^{p}  
+\lambda^{p/2} |u |^{p} ))_{Q_{\nu r}}.
\end{equation}

 \end{theorem}

{\bf Proof of Theorem \ref{theorem 06.10.6.1}}.
 We follow  the general scheme of   proving  
Theorem \ref{theorem 06.6.13.1}. We may certainly assume
that $u$ and $f$ have compact supports.
Then as in the proof of Lemma \ref{lemma 7.8.1}
we
  may assume that
$a,u,f$ are infinitely differentiable.

In that case
take a $\lambda>0$, which in the future will be sent to 0,
 take a $\zeta\in C^{\infty}_{0}(\bR^{d+1})$ such that
$\zeta=1$ on $Q_{\nu r/2}-Q_{\nu r/2}$ and $\zeta=0$ outside $Q_{\nu
r}-Q_{\nu r}$ and set
$$
g=\Div(f \zeta),\quad h=\bar{L}u-g.
$$
Next, we define $v,w^{i}$, and $\phi$ as the  unique solutions in
$W^{1,2}_{p} $ of the equations
$$
\bar{L}v -\lambda v=h,\quad
\bar{L}w^{i} -\lambda w^{i}=f^{i}\zeta,\quad \bar{L}\phi -\lambda\phi=
-\lambda u.
$$
Since $\lambda>0$, by classical theory we know that such   $v,w,\phi$ 
indeed exist, are unique, and infinitely differentiable.  

  Since $h=0$ in
$Q_{\nu r/2} $ and  
$\nu/2\geq2$, by Theorem
\ref{theorem 06.10.3.1} we obtain
$$
(|v_{x } 
-(v_{x })_{Q_{r}}|^{p})_{Q_{r}}
\leq N\nu^{-p}
(|v_{x }|^{p}   
+\lambda^{p/2} |v|^{p} )_{Q_{\nu r/2}}
$$
\begin{equation}
                                                  \label{06.10.6.2}
 \leq N\nu^{-p}
(|v_{x }|^{p}   
+\lambda^{p/2} |v|^{p} )_{Q_{\nu r }}.
\end{equation}

Furthermore,
$$
\lambda\|w\|_{L_{p}(\bR^{d+1}_{0})}
+\lambda^{1/2}\|w_{x}\|_{L_{p}(\bR^{d+1}_{0})}
$$
$$
+
\|w_{t}\|_{L_{p}(\bR^{d+1}_{0})}+\|w_{xx}\|_{L_{p}(\bR^{d+1}_{0})}
\leq N\|f\zeta\|_{L_{p}(\bR^{d+1}_{0})}.
$$
In particular, for $\psi:=\Div w$ we have
$$
\lambda^{p/2}\int_{Q_{\nu r}}|\psi|^{p}\,dxdt
+\int_{Q_{\nu r}}|\psi_{x}|^{p}\,dxdt
\leq N\int_{Q_{\nu r}}|f|^{p}\,dxdt.
$$
\begin{equation}
                                                  \label{06.10.6.5}
(|\psi_{x }|^{p}   
+\lambda^{p/2} |\psi|^{p} )_{Q_{ r}}
\leq N\nu^{d+2}(|f|^{p})_{Q_{\nu r}}.
\end{equation}

Also,
$$
\lambda\|\phi\|_{L_{p}(\bR^{d+1}_{0})}
+\lambda^{1/2}\|\phi_{x}\|_{L_{p}(\bR^{d+1}_{0})}
$$
$$
+
\|\phi_{t}\|_{L_{p}(\bR^{d+1}_{0})}+\|\phi_{xx}\|_{L_{p}(\bR^{d+1}_{0})}
\leq N\lambda\|u\|_{L_{p}(\bR^{d+1}_{0})},
$$
$$
(|\phi_{x }|^{p}   
+\lambda^{p/2} |\phi|^{p} )_{Q_{\nu r}}
\leq N\lambda^{p/2}(\nu r)^{-d-2} 
\|u\|^{p}_{L_{p}(\bR^{d+1}_{0})},
$$
\begin{equation}
                                                  \label{06.10.6.6}
(|\phi_{x }|^{p}   
+\lambda^{p/2} |\phi|^{p} )_{Q_{ r}}
\leq N\lambda^{p/2}r^{-d-2} 
\|u\|^{p}_{L_{p}(\bR^{d+1}_{0})}.
\end{equation}

Finally, we claim that
$$
u=v+\psi+\phi=:\bar{u}.
$$
Indeed, owing 
to the additional assumptions on $f$, for any multi-index $\alpha$
we have $D^{\alpha}w\in W^{1,2}_{p}$. Hence, 
$\bar{u} \in W^{1,2}_{p}$. Upon observing that
$$
\bar{L}\bar{u}-\lambda\bar{u}=
h+\Div(f\zeta)-\lambda u=\bar{L}u-\lambda u
$$
and using uniqueness we get that $\bar{u}=u$, indeed.

After that, by using
\eqref{06.10.6.2},  \eqref{06.10.6.5}, and \eqref{06.10.6.6},
 we can dominate the left-hand side of
\eqref{06.10.6.1} by a constant times
$$
\big(|v_{x }-(v_{x })_{Q_{r}}| ^{p}\big)_{Q_{r}}
+\big(|\psi_{x}| ^{p}\big)_{Q_{r}}
+\big(|\phi_{x}| ^{p}\big)_{Q_{r}}
\leq 
N\nu^{-p}
 (|v_{x }|^{p}   
+\lambda^{p/2} |v|^{p} )_{Q_{\nu r}}
$$
$$
+N\nu^{d+2} (|f|^{p})_{Q_{\nu r}}
+N\lambda^{p/2}r^{-d-2} 
\|u\|^{p}_{L_{p}(\bR^{d+1}_{0})},
$$
where
$$
 (|v_{x }|^{p}   
+\lambda^{p/2} |v|^{p} )_{Q_{\nu r}}
\leq N (|u_{x }|^{p}   
+\lambda^{p/2} |u|^{p} )_{Q_{\nu r}}
$$
$$
+N (|\psi_{x }|^{p}   
+\lambda^{p/2} |\psi|^{p} )_{Q_{\nu r}}
+N (|\phi_{x }|^{p}   
+\lambda^{p/2} |\phi|^{p} )_{Q_{\nu r}}
$$
$$
\leq 
N 
(|u_{x }|^{p}   
+\lambda^{p/2} |u|^{p} )_{Q_{\nu r}}
+N  (|f|^{p})_{Q_{\nu r}} 
+N\lambda^{p/2}  r ^{-d-2} 
\|u\|^{p}_{L_{p}(\bR^{d+1}_{0})}.
$$

Thus, the left-hand side of
\eqref{06.10.6.1} is less than
$$
N \nu^{-p}
(|u_{x }|^{p}   
+\lambda^{p/2} |u|^{p} )_{Q_{\nu r}}
+N \nu^{d+2} (|f|^{p})_{Q_{\nu r}} 
+N\lambda^{p/2}  r ^{-d-2} 
\|u\|^{p}_{L_{p}(\bR^{d+1}_{0})}
$$
and to obtain \eqref{06.10.6.1} it only remains to let
$\lambda\downarrow0$.
  The theorem is proved.

Now we can repeat what is said in \cite{Kr06}
and get the solvability of equations involving $\cL$
in $\cH^{1}_{p}$. In particular, we have the following result.
Recall that $\bR_{S}$ and $\bH^{-1}_{p}((S,T))$
are introduced in Section \ref{intro}.

\begin{theorem}
                                                  \label{theorem 7.18.2}
Let $S\in[-\infty ,\infty)$. Then there exists
  $\lambda_{0}$, depending only on
$p,d,K,\kappa$, and $\omega$,
such that,    for any
$u\in\cH^{1}_{p}(\bR^{d+1}_{S}) $ and $\lambda \geq\lambda_{0}$,   we have
\begin{equation}
                                                       \label{7.18.8}
\|u_{t}\|_{\bH^{-1}_{p}(\bR_{S}) }+\|u_{x}\|_{
L_{p}(\bR^{d+1}_{S}) }
 +\|u\|_{L_{p}(\bR^{d+1}_{S})}\leq 
N\|(\cL -\lambda)
u\|_{\bH^{-1}_{p} (\bR_{S})},
\end{equation}
where $N$ depends only on $ p,d,K,\kappa,\omega,\lambda$.
 Furthermore, for each  $\lambda\geq\lambda_{0}$ and
  $f\in \bH^{-1}_{p}(\bR_{S}) $
there is a unique
$u\in\cH^{1}_{p}(\bR_{S}) $ such that   $(\cL-\lambda) u=f$.

\end{theorem}

This is a version of Theorem 6.2 of \cite{Kr06}.
There is only one difference.   Theorem 6.2 of \cite{Kr06}
  is stated with
$\bH^{-1}_{p}$ and $L_{p}$ in place of 
$\bH^{-1}_{p}(\bR_{S})$ and $L_{p}(\bR^{d+1}_{S})$, respectively.
  Passing from the former spaces to the latter ones
is performed as in \cite{Kr06} on the basis of the fact that
the a priori estimate \eqref{7.18.8} allows one to solve
the corresponding 
equations by the method of continuity
(cf. Corollary \ref{corollary 06.10.5.1}).

\mysection{Proof of Lemma \protect\ref{lemma 06.10.9.5}}
                            \label{section 06.10.20.1} 

Although the way we proceed are similar to
what is done in \ref{section 06.10.19.1} 
the details are quite different
and the main reason for that is that
we cannot prove a natural counterpart
of Lemma \ref{lemma 06.10.5.1}.

 \begin{lemma}
                                                \label{lemma 06.10.11.2}
Let   $r\in(0,\infty)$,
$q\in(1,p]$, and assume that
\begin{equation}
                                                  \label{06.10.8.5}
\frac{1}{q}-\frac{1}{p}\leq\frac{1}{d+2}.
\end{equation}
Let $\zeta\in C^{\infty}_{0}(\bR^{d+1})$
be such that $\zeta=1$ in $Q_{r}$.
Then for any function $u$ such that
 $u\zeta\in \cH^{1}_{q}(\bR_{0}) $, 
we have $u \in L_{p}(Q_{r})$ and
$$
 \|u\|_{L_{p}(Q_{r})}\leq
N\|u\zeta\|_{\cH^{1}_{q}(\bR_{0})},
$$
where $N=N(r,d, p,q,\zeta)$.
\end{lemma}

Proof. Take a $\lambda_{0}$ which suits $\cL=D_{t}+\Delta$
in Theorem \ref{theorem 7.18.2}.
By definition  
$u\zeta=(1-\Delta)^{1/2}w$, where $w\in W^{1,2}_{q}(\bR^{d+1}_{0})$
and
$$
w_{t}+\Delta w-\lambda_{0}w=:\phi\in L_{q}(\bR^{d+1}_{0}).
$$
By applying to both parts $(1-\Delta)^{1/2}$ we see that
$$
h:=\Delta(u\zeta)+(u\zeta)_{t}-
\lambda_{0}u\zeta=(1-\Delta)^{1/2}\phi.
$$
Observe that
\begin{equation}
                                                  \label{06.10.18.5}
\|h\|_{\bH^{-1}_{q}(\bR_{0})}=
\|\phi\|_{L_{q}(\bR^{d+1}_{0})}
\leq N\|w\|_{W^{1,2}_{q}(\bR^{d+1}_{0})}
=N\|u\zeta\|_{\cH^{1}_{q}(\bR_{0})}.
\end{equation}

Next, write
$$
h=\Div f+g,\quad
g:=(1-\Delta)^{-1/2}\phi=(1-\Delta)^{-1}h,\quad f:=-g_{x}
$$
and notice that
\begin{equation}
                                                  \label{06.10.11.6}
\|f\|_{L_{q}(\bR^{d+1}_{0})}+
\|g\|_{L_{q}(\bR^{d+1}_{0})}\leq
N\|h\|_{\bH^{-1}_{q}(\bR_{0})},
\end{equation}
where $N=N(d,q)$.

Now define $v$ and $w$ as the unique solutions
from $W^{1,2}_{q}(\bR^{d+1}_{0})$
of
$$
\Delta v+v_{t}-\lambda_{0}v=g,\quad \Delta w+w_{t}-\lambda_{0}w=f.
$$
Then we have $\bar{u}:=v+\Div w\in \cH^{1}_{q}(\bR _{0})$
 since 
$$
D_{i}W^{1,2}_{q}(\bR^{d+1}_{0})
=(1-\Delta)^{1/2}[(1-\Delta)^{-1/2}D_{i}]
W^{1,2}_{q}(\bR^{d+1}_{0})
$$
$$
\subset
(1-\Delta)^{1/2}W^{1,2}_{q}(\bR^{d+1}_{0}).
$$
Furthermore,
obviously $\Delta\bar{u}+\bar{u}_{t}-\lambda_{0}\bar{u}=h$.
Since $u\zeta$ also satisfies this equation, by 
Theorem \ref{theorem 7.18.2}, we have
$\bar{u}=u\zeta$. In particular, $ u=v+\Div w$
in $Q_{r}$.

Finally, by classical results and \eqref{06.10.11.6}
and \eqref{06.10.18.5}
$$
\|v_{xx}\|_{L_{q}(\bR^{d+1}_{0})}
+\|v_{t}\|_{L_{q}(\bR^{d+1}_{0})}+\|v \|_{L_{q}(\bR^{d+1}_{0})}
\leq N\|g\|_{L_{q}(\bR^{d+1}_{0})}
\leq N\|u\zeta\|_{\cH^{1}_{q}(\bR_{0})},
$$
$$
\|w_{xx}\|_{L_{q}(\bR^{d+1}_{0})}
+\|w_{t}\|_{L_{q}(\bR^{d+1}_{0})}+\|w\|_{L_{q}(\bR^{d+1}_{0})}
\leq N\|u\zeta\|_{\cH^{1}_{q}(\bR_{0})}.
$$
It only remains to notice that,
owing to \eqref{06.10.8.5}, by Lemma II.3.3 of \cite{LSU}
$$
\|v\|_{L_{p}(\bR^{d+1}_{0})}\leq N(
\|v_{xx}\|_{L_{q}(\bR^{d+1}_{0})}
+\|v_{t}\|_{L_{q}(\bR^{d+1}_{0})}+\|v \|_{L_{q}(\bR^{d+1}_{0})}),
$$
$$
\|w_{x}\|_{L_{p}(\bR^{d+1}_{0})}\leq N(
\|w_{xx}\|_{L_{q}(\bR^{d+1}_{0})}
+\|w_{t}\|_{L_{q}(\bR^{d+1}_{0})}+\|w\|_{L_{q}(\bR^{d+1}_{0})}).
$$
The lemma is proved.

To move further to equations
in $\cH^{1}_{q,p}$ spaces we need the following
counterpart of Lemma \ref{lemma 06.10.5.1}.

 \begin{lemma}
                                                \label{lemma 06.10.11.1}
Let $r\in(0,1]$, $\nu\in(1,\infty)$,
$q\in(1,p]$, and assume \eqref{06.10.8.5}.
Let $u\in \cH^{1}_{q,loc} $, 
$f=(f^{1},...
,f^{d})$, $f^{i}, g\in L_{q,loc}$  and assume that $ \cL
u=\Div f+g$ in $Q_{\nu r}$.  
  Then $u \in L_{p}(Q_{r})$ and
\begin{equation}
                                                  \label{06.10.11.7}
r^{-1} ( |u |^{p} )^{1/p}_{Q_{r}}
\leq N  ( |f |^{q} 
+ r |g|^{q}  
+  |u_{x}|^{q}  
+ r^{-1}  |u |^{q} )^{1/q}_{Q_{\nu r}}  ,
\end{equation}
where $N=N(\nu,d,\kappa,p,q,K)$.
Furthermore, if, additionally, $f  \in L_{p,loc}$,
then $u_{x}\in L_{p}(Q_{r})$ and
\begin{equation}
                                                  \label{06.10.11.1}
 ( |u _{x}|^{p} )^{1/p}_{Q_{r}}
\leq N\big [ ( |f |^{p} )^{1/p}_{Q_{\nu r}}
+ (r|g|^{q}  
+  |u_{x}|^{q} 
+ r^{-1}( |u |^{q} )^{1/q}_{Q_{\nu r}}\big ] ,
\end{equation}
where $N=N(\nu,d,\kappa,p,q,K)$.
\end{lemma}

Proof. 
By self-similarity we may assume that $r=1$
(cf. the proof of Lemma \ref{lemma 06.10.5.1}). 
In that case
take 
$\lambda=\lambda_{0}$
which suits
   Theorem \ref{theorem 7.18.2} for both $p$
and $q$ in place of $p$ there.
Also take a $\zeta\in C^{\infty}_{0}(\bR^{d+1}_{0})$ such that
$\zeta=1$ on $Q_{1}$ and $\zeta=0$ in $\bR^{d+1}_{0}
\setminus Q_{\nu}$.
Observe that in $\bR^{d+1}_{0}$ we have
\begin{equation}
                                                  \label{06.10.11.8}
\cL(u\zeta )-\lambda u\zeta =\Div(\zeta f+\bar{f})+\bar{g},
\end{equation}
where
$$
\bar{f}^{j}= ua^{ij}\zeta_{x^{i}},\quad
\bar{g}= \zeta  g
-f^{i}\zeta_{ x^{i}}
+
u\big[\zeta_{t}+(b^{i}+\hat{b}^{i})\zeta_{ x^{i}}-\lambda \zeta\big]
+a^{ij}u_{x^{i}}\zeta_{x^{j}} .
$$
Since the $\bH^{-1}_{q}(\bR _{0})$-norm of the right-hand side of
\eqref{06.10.11.8} is less than a constant times
the right-hand side of
\eqref{06.10.11.7} we get \eqref{06.10.11.7} (with $r=1$)
by Theorem \ref{theorem 7.18.2}
 and Lemma \ref{lemma 06.10.11.2}.
Below we are going to use a trivial extension
of this result that \eqref{06.10.11.7} also holds
with $Q_{\nu'}$ in place of $Q_{1}$, where
$\nu'=(1+\nu)/2$. We will also assume that
$\zeta=0$ in $\bR^{d+1}_{0}\setminus Q_{\nu'}$.

To prove \eqref{06.10.11.1} we  want to apply 
Theorem \ref{theorem 7.18.2} again to \eqref{06.10.11.8}. By the above
the $L_{p}(\bR^{d+1}_{0})$-norm of $\zeta f+\bar{f}$
is under control. To deal with $g$
we define $v$ as the unique $W^{1,2}_{q}(\bR^{d+1}_{0})$
solution of $\Delta v+v_{t}-\lambda v=\bar{g}$.
Notice that, as in the proof of Lemma 
\ref{lemma 06.10.11.2},
\begin{equation}
                                                  \label{06.10.11.9}
\|v \|_{L_{p}(\bR^{d+1}_{0})}+
\|v_{x}\|_{L_{p}(\bR^{d+1}_{0})}
\leq N\|\bar{g}\|_{L_{q}(\bR^{d+1}_{0})}\leq NI_{1},
\end{equation}
where $I_{r}$ is the content  of the brackets
on the right in \eqref{06.10.11.1}.
Furthermore, $w=u\zeta-v$ satisfies
\begin{equation}
                                                  \label{06.10.11.10}
\cL w-\lambda w=\Div (\zeta f+\bar{f}+\hat{f})
+\hat{g},
\end{equation}
where
$$
\hat{f}^{j}=(\delta^{ij}-a^{ij})v_{x^{i}}-\hat{b}^{j}v,\quad
\hat{g}=-b^{i}v_{x^{i}}-cv.
$$
By the above the right-hand side of \eqref{06.10.11.10}
is in $\bH^{-1}_{p}(\bR_{0})$ with norm controlled by
$NI_{1}$. By Theorem \ref{theorem 7.18.2}
equation  \eqref{06.10.11.10} has a unique solution
in $\cH^{1}_{p}(\bR _{0})$ and, since
$w=u\zeta-v\in \cH^{1}_{q}(\bR _{0})$,
this certainly implies that
$w \in \cH^{1}_{p}(\bR _{0})$. Also by
Theorem \ref{theorem 7.18.2}
$$
\|w_{x}\|_{L_{p}(\bR^{d+1}_{0})}
\leq N I_{1},
$$
which along with \eqref{06.10.11.9} leads to
\eqref{06.10.11.1}. The lemma is proved.

 \begin{corollary}
                                 \label{corollary 06.10.8.3}
Let $r\in(0,1]$, $\nu\in(1,\infty)$,
$q\in(1,p]$,
$u\in \cH^{1}_{q,loc}$, 
$f=(f^{1},...
,f^{d})$, $f^{i},g\in L_{p,loc} $ and assume that $ \cL
u=\Div f+g$ in $Q_{\nu r}$.  
  Then $u,u_{x}\in L_{p}(Q_{r})$ and
$$
r^{-1} ( |u |^{p} )^{1/p}_{Q_{r}}
+( |u_{x} |^{p} )^{1/p}_{Q_{r}}
\leq N( |f |^{p} )^{1/p}_{Q_{\nu r}}
$$
\begin{equation}
                                                  \label{06.10.8.4}
+Nr( |g |^{p} )^{1/p}_{Q_{\nu r}}
+N( |u_{x} |^{q} )^{1/q}_{Q_{\nu r}}
+Nr^{-1}( |u |^{q} )^{1/q}_{Q_{\nu r}} ,
\end{equation}
where $N=N(\nu,d,\kappa,p,q,K)$.
\end{corollary}

Indeed, if our $p$ satisfies
\eqref{06.10.8.5}, then we have the result by 
Lemma \ref{lemma 06.10.11.1}. If $p$ is bigger, then we use
Lemma \ref{lemma 06.10.11.1} with $p_{1}$ in place of $p$,
where $p_{1}$ is defined by 
$q^{-1}-p_{1}^{-1}=(d+2)^{-1}$. Once we have the result
for $p_{1}$, we take $p_{1}$ as new $q$ and keep iterating as many
times as needed, each time reducing $p_{k}^{-1}$ by $(d+2)^{-1}$
until we reach first $k$ such that $p_{k}^{-1}-p^{-1}\leq(d+2)^{-1}$.

\begin{corollary}
                                           \label{corollary 06.10.9.1}
If
  $r\in(0,1]$, $q>1$, and
$u\in \cH^{1}_{q,loc}  $ 
are such that that in $Q_{2 r}$ we have
$\cL u  =0$, $b=\hat{b} =0$, and $c=0$,
then $u_{x}\in L_{p}(Q_{r})$ and
\begin{equation}
                                                  \label{06.10.9.6}
 ( |u_{x }|^{p} )^{1/p}_{Q_{r}}
\leq
N   ( |u_{x }|^{q} )^{1/q}_{Q_{2r}} ,
\end{equation}
where   
 $N $ depends only on $ d ,p,q,K,
\kappa
$, and the function $\omega$.
\end{corollary}

For $q\geq p$ equation \eqref{06.10.9.6} is obvious.
To prove it for $q\leq p$ it suffices to apply \eqref{06.10.8.4}
to $v=u-u_{Q_{\nu r}}$ in place of $u$, observe that
$\cL v=0$ in $Q_{2 r}$, and finally use
Lemma \ref{lemma 7.8.1} with $\bar{L}=\Delta+D_{t}$
for which
$$
\bar{L}u=((\delta^{ij}-a^{ij})u_{x^{i}})_{x^{j}}.
$$

{\bf Proof of Lemma \ref{lemma 06.10.9.5}}.
 We may certainly assume that
$u\in\cH^{1}_{p}$. According to Theorem 2.4 of \cite{Kr06},
 applied to the domains $ 
 (S,4)\times\bR^{d}$ for $S<4$,
on $\bR^{d+1} $ there is a function $v$ such that
it belongs to $\cH^{1 }_{p}(\bR _{S})$ for any $S$,
satisfies
\begin{equation}
                                                \label{06.10.10.1}
\cL v =\Div(fI_{Q_{\nu r}}) 
\end{equation}
in $\bR^{d+1} $, and is such that $v(t,x)=0$ for $t>4$
(observe that $\nu r\leq1$). 
 
After that we  set
$$
w=u-v
$$ 
and note for the future that
$w\in \cH^{1 }_{p,loc} $.

Again by Theorem 2.4 of \cite{Kr06} we have
$$
\int_{(0,4)\times\bR^{d}}
|v_{x }|^{p}\,dxdt\leq N\int_{Q_{\nu r}}|f|^{p}\,dxdt
$$
implying that
\begin{equation}
                                                \label{06.10.10.2}
(|v_{x }|^{p})_{Q_{\nu r}}\leq N\cA_{\nu r},\quad
(|v_{x }|^{p})_{Q_{r}}\leq N\nu^{d+2}\cA_{\nu r}.
\end{equation}

Next, since
$\cL w =0$ in $Q_{\nu r}$  and $\nu r/2\leq1$,
Corollary \ref{corollary 06.10.9.1} implies that
$w_{x}\in L_{2p}(Q_{\nu r/2})$ and
\begin{equation}
                                                \label{06.10.11.11}
 ( 
|w_{x}|^{2p} )^{1/2}_{Q_{\nu r/2}}
\leq N (
|w_{x}|^{ p})_{Q_{\nu r }}.
\end{equation}

 Upon taking an
$\bar{a}\in\bA$, setting $\bar{L}\phi=\phi_{t}
+\bar{a}^{ij}\phi_{x^{i}x^{j}}$ and noting that
$$
\bar{L}w=\Div\bar{f},\quad
\bar{f}^{j}= (\bar{a}^{ij}-a^{ij})w_{x^{i}} ,
$$
by Theorem \ref{theorem 06.10.6.1} we get
$$
\big(|w_{x }-( w_{x })_{Q_{r}}| ^{p}\big)_{Q_{r}}
\leq N\nu^{-p}
\big(|w_{ x}|^{p} \big)_{Q_{\nu r/2}} 
+N \nu^{ d+2 }
\big(  |\bar{f}|^{p}   \big)_{Q_{\nu r/2}}, 
$$
where the last term by H\"older's inequality and
\eqref{06.10.11.11} is dominated by a constant
times
$$
(|\bar{a}-a|^{2p})_{Q_{\nu r/2}}^{1/2}
(|w_{x}|^{2p})_{Q_{\nu r/2}}^{1/2}
\leq(|\bar{a}-a| )_{Q_{\nu r/2}}^{1/2}
(|w_{x}|^{ p})_{Q_{\nu r/2}} .
$$

By using an appropriate choice of $\bar{a}$ we obtain
$$
(  |\bar{f}|^{p}   \big)_{Q_{\nu r/2}}
\leq N\hat{a}^{1/2} (|w_{x}|^{ p})_{Q_{\nu r }}.
$$
 Combining the above
and observing that
$$
\big(|w_{x}|^{p}\big)_{Q_{\nu r/2}}
\leq 2^{d+2}\big(|w_{x}|^{p}\big)_{Q_{\nu r }}
\leq N\big(|u_{x}|^{p}\big)_{Q_{\nu r}}
+N\big(|v_{x}|^{p}\big)_{Q_{\nu r}}
$$
$$
\leq N\big(|u_{x}|^{p}\big)_{Q_{\nu r}}
+N\cA_{\nu r}
$$
 yield
that the left-hand side of
\eqref{06.10.9.7} is dominated by
$$
N\big(|v_{x}|^{p}\big)_{Q_{r}}
+N(|w_{x }-(w_{x })_{Q_{r} }|^{p})_{Q_{r} }
\leq N\nu^{d+2}\cA_{\nu r}
+N\nu^{-p}\big(|u_{x}|^{p}\big)_{Q_{\nu r}}
$$
$$
+N\nu^{d+2}\hat{a}^{1/2}\big(
\big(|u_{x}|^{p}\big)_{Q_{\nu r}}
+\cA_{\nu r}\big).
$$
This is almost exactly what is asserted and the lemma is proved.

\end{document}